\documentclass[11pt]{article}

\usepackage{a4,amsfonts,exscale}

\long\def\forget#1{}

\usepackage{amsmath}
\usepackage{amssymb}
\usepackage{amscd}
\usepackage{textcomp}

\usepackage{amsthm}
\theoremstyle{plain}

\theoremstyle{definition}

\theoremstyle{remark}

\usepackage{xy}
\xyoption{all}

\newdir^{ (}{{}*!/-3pt/\dir^{(}}    
\newdir^{  }{{}*!/-3pt/\dir^{}}    
\newdir_{ (}{{}*!/-3pt/\dir_{(}}    
\newdir_{  }{{}*!/-3pt/\dir_{}}

\newcounter{zahl}

\def\theenumi{(\alph{enumi})}

\def\p@enumii{\theenumi}

\newcommand{\DS}{\displaystyle}
\newcommand{\TS}{\textstyle}
\newcommand{\SC}{\scriptstyle}
\newcommand{\SSC}{\scriptscriptstyle}

\DeclareMathOperator{\End}{End}

\DeclareMathOperator{\Frob}{Frob}
\DeclareMathOperator{\Gal}{Gal}
\DeclareMathOperator{\GL}{GL}

\DeclareMathOperator{\Koh}{H}

\DeclareMathOperator{\Lie}{Lie}
\DeclareMathOperator{\Mod}{Mod}

\DeclareMathOperator{\Quot}{Frac}

\DeclareMathOperator{\Rep}{Rep}

\DeclareMathOperator{\Spm}{Sp}
\DeclareMathOperator{\Spec}{Spec}

\newcommand{\alg}{{\rm alg}}

\DeclareMathOperator{\coker}{coker}
\newcommand{\cris}{{\rm cris}}

\newcommand{\dR}{{\rm dR}}
\newcommand{\et}{{\rm\acute{e}t}}

\renewcommand{\mod}{{\rm\,mod\,}}
\DeclareMathOperator{\ord}{ord}
\newcommand{\perf}{{\rm perf}}

\newcommand{\rig}{{\rm rig}}
\DeclareMathOperator{\rk}{rk}
\newcommand{\sep}{{\rm sep}}

\renewcommand{\phi}{\varphi}
\renewcommand{\epsilon}{\varepsilon}

\usepackage{amsfonts}
\newcommand{\BOne} {{\mathchoice{\hbox{\rm1\kern-2.7pt l\kern.9pt}}
                              {\hbox{\rm1\kern-2.7pt l\kern.9pt}}
                              {\hbox{\scriptsize\rm1\kern-2.3pt l\kern.4pt}}
                              {\hbox{\scriptsize\rm1\kern-2.4pt l\kern.5pt}}}}

\newcommand{\BD}{{\mathbb{D}}}
\newcommand{\BE}{{\mathbb{E}}}
\newcommand{\BF}{{\mathbb{F}}}

\newcommand{\BH}{{\mathbb{H}}}

\newcommand{\BM}{{\mathbb{M}}}
\newcommand{\BN}{{\mathbb{N}}}

\newcommand{\BQ}{{\mathbb{Q}}}
\newcommand{\BR}{{\mathbb{R}}}

\newcommand{\BZ}{{\mathbb{Z}}}

\newcommand{\bA}{{\mathbf{A}}}
\newcommand{\bB}{{\mathbf{B}}}

\newcommand{\bD}{{\mathbf{D}}}
\newcommand{\bE}{{\mathbf{E}}}

\newcommand{\CF}{{\cal{F}}}
\newcommand{\CG}{{\cal{G}}}
\newcommand{\CH}{{\cal{H}}}

\newcommand{\CO}{{\cal{O}}}

\newcommand{\CR}{{\cal{R}}}

\newcommand{\FM}{{\mathfrak{M}}}

\newcommand{\FS}{{\mathfrak{S}}}

\newcommand{\Fm}{{\mathfrak{m}}}

\newcommand{\Fp}{{\mathfrak{p}}}
\newcommand{\Fq}{{\mathfrak{q}}}

\newcommand{\shortonto}{\mbox{\mathsurround=0pt \;$\to \hspace{-0.8em} \to$\;}}

\newcommand{\es}{\enspace}
\newcommand{\open}{^\circ}

\newcommand{\mal}{^{\SSC\times}}

\newcommand{\ul}[1]{{\underline{#1}}}

\newcommand{\wh}[1]{{\widehat{#1}}}
\newcommand{\wt}[1]{{\widetilde{#1}}}

\usepackage{ifthen}

\newcommand{\invlim}[1][]{\ifthenelse{\equal{#1}{}}
{{\DS \lim_{\longleftarrow}\;}}
{{\DS \lim_{\underset{#1}{\longleftarrow}}\;}}
}

\newcommand{\dirlim}[1][]{\ifthenelse{\equal{#1}{}}
{{\DS \lim_{\longrightarrow}\;}}
{{\DS \lim_{\underset{#1}{\longrightarrow}}\;}}
}

\newcommand{\dbl}{{\mathchoice{\mbox{\rm [\hspace{-0.15em}[}}
                              {\mbox{\rm [\hspace{-0.15em}[}}
                              {\mbox{\scriptsize\rm [\hspace{-0.15em}[}}
                              {\mbox{\tiny\rm [\hspace{-0.15em}[}}}}
\newcommand{\dbr}{{\mathchoice{\mbox{\rm ]\hspace{-0.15em}]}}
                              {\mbox{\rm ]\hspace{-0.15em}]}}
                              {\mbox{\scriptsize\rm ]\hspace{-0.15em}]}}
                              {\mbox{\tiny\rm ]\hspace{-0.15em}]}}}}
\newcommand{\dpl}{{\mathchoice{\mbox{\rm (\hspace{-0.15em}(}}
                              {\mbox{\rm (\hspace{-0.15em}(}}
                              {\mbox{\scriptsize\rm (\hspace{-0.15em}(}}
                              {\mbox{\tiny\rm (\hspace{-0.15em}(}}}}
\newcommand{\dpr}{{\mathchoice{\mbox{\rm )\hspace{-0.15em})}}
                              {\mbox{\rm )\hspace{-0.15em})}}
                              {\mbox{\scriptsize\rm )\hspace{-0.15em})}}
                              {\mbox{\tiny\rm )\hspace{-0.15em})}}}}

\newcommand{\BFZ}{{\BF_q\dpl\zeta\dpr}}
\newcommand{\con}[1][]{{\mathchoice
           {\TS\langle\frac{z}{\zeta^{#1}}\rangle[\frac{1}{z}]}
           {\TS\langle\frac{z}{\zeta^{#1}}\rangle[\frac{1}{z}]}
           {\SC\langle\frac{z}{\zeta^{#1}}\rangle[\frac{1}{z}]}
           {\SSC\langle\frac{z}{\zeta^{#1}}\rangle[\frac{1}{z}]}}}
\newcommand{\ancon}[1][]{{\mathchoice
           {\TS\langle\frac{z}{\zeta^{#1}},\frac{1}{z}\}}
           {\TS\langle\frac{z}{\zeta^{#1}},\frac{1}{z}\}}
           {\SC\langle\frac{z}{\zeta^{#1}},\frac{1}{z}\}}
           {\SSC\langle\frac{z}{\zeta^{#1}},\frac{1}{z}\}}}}

\newcommand{\dotBD}{\vbox{\hbox{\kern2pt\bf.}\vskip-4.5pt\hbox{$\BD$}}}

\def\?{\ 

???\ \immediate\write16{}

\immediate\write16{Warning: There was still a question mark . . . }
\immediate\write16{}}

\def\longto{\longrightarrow}

\def\isoto{\stackrel{\;\sim\es}{\longrightarrow}}

\newbox\mybox
\def\arrover#1{\mathrel{
       \setbox\mybox=\hbox spread 1.4em{\hfil$\scriptstyle#1$\hfil}
       \vbox{\offinterlineskip\copy\mybox
             \hbox to\wd\mybox{\rightarrowfill}}}}

\newlength{\columnMixed}
\newlength{\columnEqual}
\newcommand{\comp}[2]{\noindent\begin{tabular}{l|r} \mbox{ } \\[-5.9mm]\mbox{ }\\[-3mm] \parbox[t]{\columnMixed}{#1} & \parbox[t]{\columnEqual}{#2} \\[-3mm] \mbox{ } \\[-5.9mm]\mbox{ }\end{tabular}}

\newcommand{\commt}[1]{\noindent\begin{tabular}{p{0.97\textwidth}} #1 \end{tabular} }

\begin{document}


\author{Urs Hartl
\footnote{The author acknowledges support of the Deutsche Forschungsgemeinschaft in form of DFG-grant HA3006/2-1}
}

\title{A Dictionary between Fontaine-Theory \\ and its Analogue in Equal Characteristic}

\maketitle

\begin{abstract}
\noindent
In this survey we explain the main ingredients and results of the analogue of Fontaine-Theory in equal positive characteristic which was recently developed by Genestier-Lafforgue and the author.

\end{abstract}


\bigskip

%
%

\section*{Introduction}
\addcontentsline{toc}{section}{Introduction}

The aim of Fontaine-Theory is to classify $p$-adic Galois representations of $p$-adic fields and to attach various invariants to them. To achieve this, Fontaine constructed for a $p$-adic field $K$ several $\BQ_p$-algebras with $\CG_K:=\Gal(K^\alg/K)$-action and additional structure. Then for a representation $\CG_K\to\GL(V)$ in a finite dimensional $\BQ_p$-vector space $V$ and any such $\BQ_p$-algebra $B$ the inherited additional structure on $(B\otimes_{\BQ_p}V)^{\CG_K}$ provides these invariants. This approach was extremely successful. In its application to geometry it allowed to recover all the cohomological invariants attached to a smooth proper variety $X$ over $K$ solely from the \'etale cohomology $\Koh^\bullet_\et(X\times_KK^\alg,\BQ_p)$ of $X$ (which is a $\CG_K$-representation). After contributions by Grothendieck, Tate, Fontaine, Lafaille, Messing, an many others the later was finally accomplished by Faltings~\cite{Faltings1,Faltings} and Tsuji~\cite{Tsuji}. See \cite{Illusie} for a survey.

Inspired by the close parallel between number fields and function fields, an analogue for Fontaine-Theory in equal characteristic was developed by Genestier-Lafforgue~\cite{GL} and the author~\cite{HartlPSp}. There the local field $\BQ_p$ is replaced by the Laurent series field $\BF_q\dpl z\dpr$ over the finite field $\BF_q$. But the theory of Galois representations $\CG_L:=\Gal(L^\sep/L)\to\GL_n\bigl(\BF_q\dpl z\dpr\bigr)$ for finite extensions $L$ of $\BF_q\dpl z\dpr$ is spoiled by the following two facts. Firstly by the Ax-Sen-Tate Theorem the fixed field of $\CG_L$ inside the completion $C$ of an algebraic closure of $L$ is much larger than $L$; see \ref{SectAx} below. And secondly all the Tate twists $C(n)$ are isomorphic as $\CG_L$-modules as was observed by Anderson~\cite{Anderson2}; see \ref{SectTateTwists}. Therefore it was proposed in \cite{GL,HartlPSp} to discard $\CG_L$-representations and to view in equal characteristic so called \emph{local shtuka} as the appropriate analogue of (crystalline) Galois representations; see \ref{SectRigLocalShtuka}.

This is consistent with the following facts. The link between crystalline Galois representations and other structures like filtered isocrystals first appeared via the mediator $p$-divisible groups. Fontaine-Theory then provided the direct link. Now the category of local shtuka over a formal $\BF_q\dbl z\dbr$-scheme $S$ is anti-equivalent to the category of \emph{divisible local Anderson modules} over $S$, see \ref{SectBTGroups}, and the analogues of a special class of these modules in mixed characteristic are $p$-divisible groups. Also with a divisible local Anderson module or a local shtuka over $\CO_L$ one can associate a Galois representation \mbox{$\CG_L\to\GL_n\bigl(\BF_q\dpl z\dpr\bigr)$} like one associates with a $p$-divisible group over $\CO_K$ the crystalline Galois representation on its rational Tate module; see \ref{SectTateModules}. So a priory it seems natural to use local shtuka as a mediator between Galois representations $\CG_L\to\GL_n\bigl(\BF_q\dpl z\dpr\bigr)$ and other structures. But then it turns out that the notion of local shtuka is so general (in contrast to the notion of $p$-divisible group) that one can dispense at all with Galois representations in mixed characteristic and only work with local shtuka. Striking evidence for this statement is given in \ref{SectWA=>A} and \ref{SectRZConj}.

\smallskip

This article contains no new results. It is rather a survey of the analogue of Fontaine-Theory in equal characteristic~\cite{GL,HartlPSp}. It has various purposes. First of all, a dictionary between the arithmetic of number fields and function fields was begun by Goss~\cite{GossDict}. This article is meant to be a sequel to Goss' dictionary. Secondly for those readers familiar with Fontaine-Theory it should serve as an introduction to the equal characteristic theory which is somewhat simpler but not less fascinating. Thirdly we reveal a hidden geometric interpretation of the rings in Fontaine-Theory as functions on the ``unit disc with coordinate $p$''. This interpretation is inspired  by the equal characteristic counterparts of these rings. And finally we hope that this interpretation might serve as an Ariadne's thread for those who want to learn Fontaine-Theory and experience all its rings as a kind of maze.

As with most dictionaries the reader will hardly want to read it from first to last page. We suggest to look at Section~\ref{ChaptFontainesRings}, which contains the definitions of the various rings of Fontaine-Theory and their analogues, only when these rings are needed. Instead the reader should focus on the remaining sections in which we explain the analogue of Fontaine-Theory in equal characteristic.

\bigskip

\bigskip

\noindent
\begin{tabular}[t]{p{\columnMixed}|p{\columnEqual}}
\bfseries\Large Mixed Characteristic & \hfill\bfseries\Large Equal Characteristic \\ 
\hline
\end{tabular}

\nopagebreak
\begin{center}

\section{The Arithmetic Ground Field} \label{Chapter1}

\comp{$\BQ_p$ the $p$-adic numbers and \\
$\BZ_p$ the $p$-adic integers}
{$\BF_q\dpl z\dpr$ the Laurent series field and \\
$\BF_q\dbl z\dbr$ the power series ring in the variable $z$ over the finite field $\BF_q$ with $q$ elements}

\comp{These rings carry the $p$-adic topology.}
{These rings are equipped with the $z$-adic topology.}

\subsection{Witt Vectors} \label{SectWitt}

\comp{The functor $R\mapsto W(R)$ which assigns to a perfect $\BF_p$-algebra $R$ the ring of $p$-typical \emph{Witt vectors} over $R$, see \cite{Harder,Serre},}
{corresponds to the functor $R\mapsto R\dbl z\dbr$ which assigns to an $\BF_q$-algebra $R$ the power series ring over $R$ in the variable $z$.}

\comp{The \emph{Frobenius lift} $\phi=W(\Frob_p)\in\End\bigl(W(R)\bigr)$ corresponds to}
{$\sigma:R\dbl z\dbr\to R\dbl z\dbr\,,\;\sum_i b_iz^i\mapsto \sum_i b_i^q z^i$\\[1mm]
(where $b_i\in R$), also called the \emph{Frobenius lift}.}

\comp{For an element $b\in R$ let $[b]\in W(R)$ be the \emph{Teich\-m\"uller representative} of $b$, see \cite{Harder,Serre}. The map $b\mapsto[b]$ is multiplicative but \emph{not} compatible with addition.}
{For an element $b\in R$ the element $b\cdot z^0\in R\dbl z\dbr$ is called the \emph{Teichm\"uller representative} of $b$. The map $b\mapsto b\cdot z^0$ is a ring homomorphism.}

 \comp{If $\Spec R$ is connected $W(R)^{\phi=1}=\BZ_p$.}
{If $\Spec R$ is connected $R\dbl z\dbr^{\sigma=1}=\BF_q\dbl z\dbr$.}

\subsection{The Two Roles of $p$} \label{SectTwoRolesOfP}

\comp{The number $p$ enters in Fontaine-Theory in a twofold way
  \begin{itemize}\itemsep0ex
  \item as uniformizing parameter of $\BZ_p$ and
  \item as element of the base rings or fields (which are $\BZ_p$-algebras) over which the arithmetic objects like $p$-divisible groups or Galois representations are defined.
  \end{itemize}}
{The necessity to separate the two roles of $p$ in equal characteristic and to work with two indeterminants $z$ and $\zeta$ was first pointed out by Anderson~\cite{Anderson}. So we let
  \begin{itemize}
  \item $z$ be the uniformizing parameter of $\BF_q\dbl z\dbr$,
  \item $\zeta$ be the element of the base rings (which are $\BF_q\dbl z\dbr$-algebras, $\zeta$ is the image of $z$).
  \end{itemize}}

\commt{Note that the natural number $p$ never can act on a module or vector space as anything else than the scalar $p$, whereas there is no such restriction on $z$. This makes the distinction between $z$ and $\zeta$ possible. Strictly speaking a distinction between the two roles of $p$ is also searched for in Fontaine-Theory where an object, called $[\tilde p]$ by Colmez~\cite{Colmez}, is constructed that behaves like $\zeta$, whereas $p$ behaves like $z$; see \ref{SectFieldOfNorms} and \ref{SectBdR}.}

\subsection{The Cyclotomic Character} \label{SectCyclo}

\comp{\addtolength{\baselineskip}{.7mm}For $n\in\BN_0$ let $\epsilon^{(n)}\in\BQ_p^\alg$ be a primitive $p^n$-th root of unity with $(\epsilon^{(n+1)})^p=\epsilon^{(n)}$. Let $\BQ_{p,\infty}:=\BQ_p(\epsilon^{(n)}:n\in\BN)$. Then there is an isomorphism of topological groups
\[
\chi:\Gal\bigl(\BQ_{p,\infty}/\BQ_p\bigr) \isoto \BZ_p\mal
\]
with $\gamma(\epsilon^{(n)})=(\epsilon^{(n)})^{\chi(\gamma)\mod p^n}$ for $\gamma$ in the Galois group. $\chi$ is called the \emph{cyclotomic character}.}
{\addtolength{\baselineskip}{.7mm} For $n\in\BN_0$ let $t_n\in\BF_q\dpl\zeta\dpr^\sep$ be solutions of the equations $t_0^{q-1}=-\zeta$ and $t_n^q+\zeta t_n=t_{n-1}$. Let $\BF_q\dpl\zeta\dpr_\infty:=\BF_q\dpl\zeta\dpr(t_n:n\in\BN_0)$ and define $t_{\SSC +}:=\sum_{n=0}^\infty t_n z^n\in\BF_q\dpl\zeta\dpr_\infty\dbl z\dbr$. Then there is an isomorphism of topological groups
\[
\chi:\Gal\bigl(\BF_q\dpl\zeta\dpr_\infty/\BF_q\dpl\zeta\dpr\bigr)\isoto\BF_q\dbl z\dbr\mal
\]
satisfying $\gamma(t_{\SSC +}):=\sum_{n=0}^\infty\gamma(t_n)z^n=\chi(\gamma)\cdot t_{\SSC +}$ in $\BFZ_\infty\dbl z\dbr$ for $\gamma$ in the Galois group. In view of \ref{SectLubinTate} below, $\chi$ is called the \emph{cyclotomic character}. (The existence of $\chi$ follows from the fact that by construction $\sigma(t_{\SSC +})=(z-\zeta)\cdot t_{\SSC +}$. Hence $\chi(\gamma):=\frac{\gamma(t_{\SSC +})}{t_{\SSC +}}$ is $\sigma$-invariant, that is $\chi(\gamma)\in\BF_q\dbl z\dbr\mal$. Furthermore, $\chi$ is an isomorphism because $t_{n-1}$ is a uniformizing parameter of $\BF_q\dpl\zeta\dpr(t_0,\ldots,t_{n-1})$ and so the equations defining the $t_n$ are irreducible by Eisenstein.)}

\subsection{The Ax-Sen-Tate Theorem} \label{SectAx}

\commt{Let $K$ be a field which is complete with respect to a valuation of rank one and let $C$ be the completion of an algebraic closure of $K$. By continuity the absolute Galois group $\CG_K:=\Gal(K^\sep/K)$ of $K$ acts on $C$. The Ax-Sen-Tate Theorem~\cite{Ax} states that the fixed field in $C$ of this action equals the completion of the perfection of $K$.}

\comp{If $K$ is an extension of $\BQ_p$ the fixed field is $C^{\CG_K}=K$.}
{If $K$ is an extension of $\BF_q\dpl\zeta\dpr$ the fixed field is much larger than $K$ unless $K$ is already perfect.}

\subsection{The Tate Twists} \label{SectTateTwists}

\comp{Let $K_\infty:=K(\epsilon^{(n)}:n\in\BN)$ in the notation of \ref{SectCyclo} and \ref{SectAx}. Then $\Gal(K_\infty/K)$ is a quotient of $\CG_K$ and a subgroup of $\Gal(\BQ_{p,\infty}/\BQ_p)$.}
{Let $K_\infty:=K(t_n:n\in\BN_0)$ in the notation of \ref{SectCyclo} and \ref{SectAx}. Then $\Gal(K_\infty/K)$ is a quotient of $\CG_K$ and a subgroup of $\Gal\bigl(\BFZ_\infty/\BFZ\bigr)$.}

\comp{For $m\in\BZ$ let $C(m)$ be the field $C$ on which $\gamma\in \CG_K$ acts via
\[
x\es\mapsto\es\chi(\gamma)^m\cdot \gamma(x) \quad\text{for }x\in C\,.
\]
The $C(m)$ are called the \emph{Tate twists} of $C$.}
{For $m\in\BZ$ let $C(m)$ be the field $C$ on which $\gamma\in \CG_K$ acts via
\[
x\es\mapsto\es\bigl(\chi(\gamma)|_{z=\zeta}\bigr)^m\cdot \gamma(x) \quad\text{for }x\in C\,.
\]
The $C(m)$ are called the \emph{Tate twists} of $C$.}

\comp{If $K$ is discretely valued the Tate twists are mutually non-isomorphic by \cite[Theorem 2, p.\ 176]{Tate2}.}
{Anderson~\cite{Anderson2} observed that 
all Tate twists are isomorphic! Namely, multiplication with $\bigl(t_{\SSC +}|_{z=\zeta}\bigr)^{m}=\bigl(\sum_{n\ge0}t_n\zeta^n)^{m}\in C$ defines an isomorphism $C(m)\isoto C(0)$ of $\CG_K$-modules.}

\section{Fontaine's Rings} \label{ChaptFontainesRings}

\commt{We describe most of the rings encountered in Fontaine Theory and their equal characteristic analogues. Many of these analogues are Laurent series rings in the variable $z$ over complete extensions of $\BFZ$ with varying convergence conditions. Therefore we view them as rings of (rigid) analytic functions on suitable subsets of the unit disc with coordinate $z$. We want to advertise the point of view that the corresponding rings from Fontaine Theory have a geometric interpretation as rings of analytic functions in the ``variable $p$'' on certain subsets of the ``unit disc with coordinate $p$''; see \ref{SectNonPerfect} -- \ref{SectBdR}.}

\comp{We use the notation of Colmez~\cite{Colmez} where proofs and further references can be found. Let $K$ be a \emph{complete discretely valued field extension of $\BQ_p$ with perfect residue field} $k$ and let $C$ be the completion of an algebraic closure of $K$. With the notation of \ref{SectCyclo} let $K_n:=K(\epsilon^{(n)})$ and $K_\infty:=K(\epsilon^{(n)}:n\in\BN_0)$.}
{We use the notation from \cite{HartlPSp}.  Let $L$ be a field extension of $\BFZ$ which is complete with respect to an absolute value $|\,.\,|:L\to\BR_{\geq0}$ extending the absolute value on $\BFZ$. There is no assumption on the residue field $\ell$ or on the value group of $L$. Let $C$ be the completion of an algebraic closure of $L$ and let $R$ be an \emph{affinoid $L$-algebra} with $L$-Banach norm $|\,.\,|$\,; see \cite{Bosch,BGR}. For instance the absolute case $R=L$ is allowed.}

\subsection{The Field of Norms} \label{SectFieldOfNorms}
\nopagebreak

\comp{One defines the following rings}{and their equal characteristic analogues.}

\comp{$\wt\bE^+:=\bigl\{\,x=(x^{(n)})_{n\in\BN_0}: x^{(n)}\in\CO_C,$\\
\mbox{ }\hfill $(x^{(n+1)})^p=x^{(n)}\,\bigr\}$\qquad\mbox{ }\\
$\wt\bE:=\Quot(\wt\bE^+)$ is an algebraically closed field of characteristic $p$ complete with respect to the valuation\\
$v_\bE:\wt\bE\to\BR\cup\{\infty\}\,,\,v_\bE(x):=v_p(x^{(0)})$}
{$\CO_C\cong\bigl\{\,x=(x^{(n)})_{n\in\BN_0}: x^{(n)}\in\CO_C,$\\
\mbox{ }\hfill $(x^{(n+1)})^q=x^{(n)}\,\bigr\}$\qquad\mbox{ }\\
$C=\Quot(\CO_C)$\\
$v_C:C\to\BR\cup\{\infty\}$ the valuation with $v_C(\zeta)=1$. Since we want to extend the theory to arbitrary affinoid $L$-algebras $R$ we prefer to work with the associated absolute value \\
$|\,.\,|:C\to\BR_{\ge0}\,,\,|x|=|\zeta|^{v_C(x)}$ on $C$.}

\comp{The absolute Galois group $\CG_K$ of $K$ acts component-wise on $\wt\bE^+$ and hence on $\wt\bE$.}
{The absolute Galois group $\CG_L$ of $L$ acts on $C$ but this action is not well behaved; see \ref{SectAx}, \ref{SectTateTwists}, and Section~\ref{ChaptGaloisRep}.}

\comp{$\bE^+_K:=\{\,x\in\wt\bE^+:\text{for }n\gg0 \text{ there exists }$\\
\mbox{ }\hfill $\hat x_n\in\CO_{K_n} \text{ with } v_\bE(x^{(n)}-\hat x_n)\ge 1\,\}$\quad\mbox{ }\\[1mm]
$\bE_K:=\Quot(\bE^+_K)$ is called the \emph{Field of Norms} of $K$. It is a complete discretely valued field of characteristic $p$ with perfect residue field. }
{Since the analogue of $\wt\bE$ is $C$ we may take as the analogue of $\bE_K$ the field $L$ in the absolute case or even the ring $R$ in the relative case.}

\comp{The \emph{Theory of the field of norms} states that
\[
\CH_K\es:=\es\Gal(K^\sep/K_\infty)\es=\es\Gal(\bE_K^\sep/\bE_K)\,.
\]
The reader should note that there is also a relative version of this theory by Andreatta and Brinon~\cite{Andreatta,Brinon,AB} for \emph{certain} affinoid $K$-algebras.}
{In equal characteristic $L$ (or even $R$) appears as the analogue of both $\bE_K$ and $K_\infty$. So the Theory of the field of norms is trivial in the absolute case and in the relative case for \emph{arbitrary} affinoid $L$-algebras. The analogue of $\CH_K$ is $\CH_L:=\CG_L$.}

\comp{The Theory of the field of norms is fundamental to Fontaine's theory of $p$-adic Galois representations of $K$. It allows to break up  a representation of $\CG_K$ into a representation of $\CH_K$ which corresponds to a $\phi$-module over $\bB_K$, see \ref{SectPhiGamma}, on which the remaining piece $\Gamma_K:=\CG_K/\CH_K\cong\Gal(K_\infty/K)$ of $\CG_K$ acts. On the other hand this splitting up along $K_\infty$ is unavoidable since it is impossible to identify the whole Galois group $\CG_K$ with the absolute Galois group of an appropriate field of characteristic $p$.}
{Strictly speaking $L_\infty=L(t_n:n\in\BN_0)$ is the proper analogue of both $\bE_K$ and $K_\infty$ and the Theory of the field of norms is trivial for $L_\infty$. But since this is even true for $L$ there is no need to split up Galois representations along $L_\infty$. Moreover, the notion of $z$-adic Galois representation of $L$ is not well behaved; see Section~\ref{ChaptGaloisRep}. Therefore it is preferable to work with different structures (as local shtuka, see \ref{SectLocalShtuka}) from the start. And for those structures again it is possible and better to define them over $L$ instead of $L_\infty$. For these reasons we view $L$ as the analogue of $\bE_K$.}

\comp{$\wt\bE_K$ the completion of the perfection of $\bE_K$ also equals the fixed field of $\CH_K$ in $\wt\bE$.}{$L^p:=\wh{L^\perf}$ the completion of the perfection of $L$ also equals the fixed field of $\CG_L$ in $C$.}

\comp{$\wt\bE^+_K$ the ring of integers of $\wt\bE_K$.}
{$\CO_{L^p}$ the ring of integers of $L^p$.}

\comp{Fix an element $\tilde p\in\wt\bE^+$ with $\tilde p^{(0)}=p$.\\[1mm]
Let $\epsilon:=(\epsilon^{(n)})_n\in\bE^+_{\BQ_p}$ be the sequence defined in \ref{SectCyclo} and set $\bar\pi:=\epsilon-1$.}{The analogue of $\tilde p$ is $\zeta\in L$.}

\subsection{The Witt and Power Series Rings} \label{SectWittPower}

\comp{Consider the rings}{Their equal characteristic analogues are}

\comp{$\wt\bA^+:=W(\wt\bE^+)$ the ring of Witt vectors,}
{$\CO_C\dbl z\dbr$ the ring of power series in $z$,}

\comp{$\wt\bB^+:=\wt\bA^+[\frac{1}{p}]$,}{$\CO_C\dbl z\dbr[\frac{1}{z}]$,}

\comp{$\wt\bA:=W(\wt\bE)$,}{$C\dbl z\dbr$,}

\comp{$\wt\bB:=\wt\bA[\frac{1}{p}]$ the fraction field of $\wt\bA$.}{$C\dpl z\dpr=C\dbl z\dbr[\frac{1}{z}]$ the fraction field of $C\dbl z\dbr$.}

\comp{For example $\pi:=[\epsilon]-1\in\wt\bA^+$ where $[\,.\,]$ denotes the Teichm\"uller representative, see \ref{SectWitt}.}{\mbox{ }}

\comp{There are functions $f_i:\wt\bB\to\wt\bE$ defined by the equation $\DS x=\sum_{i\gg-\infty}^\infty p^i[f_i(x)]$ in $\wt\bB$.}
{There are functions $f_i:C\dpl z\dpr\to C$ defined by the equation $\DS x=\sum_{i\gg-\infty}^\infty f_i(x)\,z^i$ in $C\dpl z\dpr$.}


\comp{$\wt\bA_K^+:=W(\wt\bE_K^+)$ the ring of Witt vectors,}
{$\CO_{L^p}\dbl z\dbr$ the ring of power series in $z$,}

\comp{$\wt\bB_K^+:=\wt\bA_K^+[\frac{1}{p}]$,}{$\CO_{L^p}\dbl z\dbr[\frac{1}{z}]$,}

\comp{$\wt\bA_K:=W(\wt\bE_K)$,}{$L^p\dbl z\dbr$,}

\comp{$\wt\bB_K:=\wt\bA_K[\frac{1}{p}]$ the fraction field of $\wt\bA_K$.}{$L^p\dpl z\dpr$ the fraction field of $L^p\dbl z\dbr$.}

\comp{By Witt vector functoriality these rings inherit a $\CG_K$-action commuting with the Frobenius lift $\phi$.}
{By functoriality these rings inherit a $\CG_L$-action which is also flawed; see Section~\ref{ChaptGaloisRep}.}

\subsection{The Non-Perfect Rings} \label{SectNonPerfect}

\comp{Set ${K_0}:=W(k)[\frac{1}{p}]$ and define \\[1mm]
$\DS \bA_{K_0}:=\bigl\{\,\sum_{n\in\BZ}a_n\pi^n:\;a_n\in W(k),\,v_p(a_n)\to\infty$\\[-2mm]
\mbox{ }\hfill for $n\to+\infty\,\bigr\}$\,,\mbox{ }\\
$\bB_{K_0}:=\bA_{K_0}[\frac{1}{p}]$ is naturally a subfield of $\wt\bB$ with residue field $\bE_{K_0}$.}{\mbox{ }}

\comp{Let $\bB_K$ be the unique finite extension of $\bB_{K_0}$ contained in $\wt\bB$ whose residue field is $\bE_K$ and let}{Its analogues are $L\dpl z\dpr$ or even $R\dbl z\dbr[\frac{1}{z}]$ in the relative situation,}

\comp{$\bA_K$ be the ring of integers of $\bB_K$ and}{respectively$L\dbl z\dbr$ or $R\dbl z\dbr$.}

\comp{$\bB^+_K:=\wt\bB^+\cap\bB_K$,}{$\CO_L\dbl z\dbr[\frac{1}{z}]$ or $R\open\dbl z\dbr[\frac{1}{z}]$ where $R\open$ is an admissible formal $\CO_L$-algebra in the sense of Raynaud~\cite{FRG} with $R\open\otimes_{\CO_L}L\cong R$. When considering $R\open$ we always assume that the $L$-Banach norm $|\,.\,|$ on $R$ satisfies $R\open=\{b\in R:|b|\le1\}$.}

\comp{$\bA^+_K:=\wt\bA^+\cap\bA_K$.}{$\CO_L\dbl z\dbr$ or $R\open\dbl z\dbr$. The elements of these rings converge on the whole (relative) open unit disc $\{|z|<1\}$ and are bounded by $1$.}

\comp{Let $\bB\subset\wt\bB$ be the completion of the maximal unramified extension of $\bB_K$ and}
{Their analogues are $L^\sep\dpl z\dpr$ and}

\comp{$\bA:=\wt\bA\cap\bB$.}{$L^\sep\dbl z\dbr$. }

\comp{By continuity $\CH_K=\Gal(\bE_K^\sep/\bE_K)$ acts on $\bB$ with $\bB^{\CH_K}=\bB_K$.}
{By continuity $\CG_L=\Gal(L^\sep/L)$ acts on $L^\sep\dpl z\dpr$ with $L^\sep\dpl z\dpr^{\CG_L}=L\dpl z\dpr$.}

\comp{The $\phi$-invariants of these rings are \\[1mm]
$(\bA^+_K)^{\phi=1}=(\bA_K)^{\phi=1}=(\wt\bA^+)^{\phi=1}=\wt\bA^{\phi=1}=\BZ_p$,}
{The $\sigma$-invariants of these rings are \\[1mm]
$R\open\dbl z\dbr^{\sigma=1}=R\dbl z\dbr^{\sigma=1}=\BF_q\dbl z\dbr$ if $\Spec R$ is connected and}

\comp{$(\bB^+_K)^{\phi=1}=(\bB_K)^{\phi=1}=(\wt\bB^+)^{\phi=1}=\wt\bB^{\phi=1}=\BQ_p$.}
{$R\open\dbl z\dbr[\frac{1}{z}]^{\sigma=1}=R\dbl z\dbr[\frac{1}{z}]^{\sigma=1}=\BF_q\dpl z\dpr$ if $\Spec R$ is connected (in particular for $R=L$ or $L^\sep$ or $C$).}

\subsection{The Overconvergent Rings} \label{SectOverconRings}

\comp{Let $\bar r\in\BR_{>0}$, set $r=1/\bar r$ and define}{Let $r\in\BR_{>0}$, set $\bar r=1/r$ and define}

\comp{$\DS\wt\bA^{(0,\bar r]}:=\bigl\{\,x\in\wt\bA:\lim_{i\to+\infty}\!v_\bE(f_i(x))+i/\bar r=\infty\,\bigr\}$,}{$C\langle\frac{z}{\zeta^r}\rangle:=\DS\bigl\{\,\sum_{i=0}^\infty b_i z^i\in C\dbl z\dbr: \lim_{i\to+\infty}\!|b_i|\,|\zeta|^{ri} = 0\,\bigr\}$.\\
Its elements converge on the disc $\{|z|\le|\zeta|^r\}$.}

\comp{$\wt\bB^{(0,\bar r]}:=\wt\bA^{(0,\bar r]}[\frac{1}{p}]$, and }{In $C\langle\frac{z}{\zeta^r}\rangle[\frac{1}{z}]$ we allow a pole at $z=0$.}

\comp{their $K$-rational versions \\[1mm]
$\bA^{(0,\bar r]}_K:=\bA_K\cap\wt\bA^{(0,\bar r]}$ and \\[1mm]
$\bB^{(0,\bar r]}_K:=\bB_K\cap\wt\bB^{(0,\bar r]}=\bA^{(0,\bar r]}_K[\frac{1}{p}]$.}
{The $L$- (or $R$)-rational versions are \\[1mm]
$L\langle\frac{z}{\zeta^r}\rangle$ or even $R\langle\frac{z}{\zeta^r}\rangle$, and \\[1mm]
$L\langle\frac{z}{\zeta^r}\rangle[\frac{1}{z}]$ or even $R\langle\frac{z}{\zeta^r}\rangle[\frac{1}{z}]$, which are defined as the preceding rings except that one must take $|\,.\,|$ to be the $L$-Banach norm on $R$ fixed in the beginning. (Since the topology of $R$ does not depend on the choice of the norm also $R\langle\frac{z}{\zeta^r}\rangle$ and $R\con[r]$ do not depend on this choice.)}

\comp{The rings $\wt\bA^{(0,\bar r]}$, $\bA_K^{(0,\bar r]}$, $\wt\bB^{(0,\bar r]}$ and $\bB_K^{(0,\bar r]}$ are principal ideal domains.}
{The rings $C\langle\frac{z}{\zeta^r}\rangle$, $L\langle\frac{z}{\zeta^r}\rangle$, $C\con[r]$ and $L\con[r]$ are principal ideal domains.}

\comp{There is a natural valuation \\[1mm]
$v^{(0,\bar r]}(x):=\inf\bigl\{\,v_\bE(f_i(x))+i/\bar r: i\in\BN_0\,\bigr\}$\\[1mm]
 on $\wt\bB^{(0,\bar r]}$ with respect to which $\wt\bA^{(0,\bar r]}$ and $\bA_K^{(0,\bar r]}$ are complete. We also consider}
{There is a natural norm \\[1mm]
$\|\sum_{i\in\BN_0}b_iz^i\|_r:=\sup\bigl\{\,|b_i|\,|\zeta|^{ri}:i\in\BN_0\,\bigr\}$\\[1mm]
on $R\langle\frac{z}{\zeta^r}\rangle$ with respect to which it is complete. We also consider}

\comp{$\DS\wt\bA^\dagger:=\bigcup_{\bar r\to0} \wt\bA^{(0,\bar r]}$,}{$\DS \bigcup_{r\to\infty}\TS C\langle\frac{z}{\zeta^r}\rangle$,}

\comp{$\DS\wt\bB^\dagger:=\bigcup_{\bar r\to0} \wt\bB^{(0,\bar r]}$ the subfield of $\wt\bB$ of \emph{overconvergent elements},}{$\DS \bigcup_{r\to\infty}\TS C\langle\frac{z}{\zeta^r}\rangle[\frac{1}{z}]$ the subfield of $C\dpl z\dpr$ of \emph{overconvergent elements} which converge on some disc with small enough radius and have a pole at $z=0$,}

\comp{$\DS\bA_K^\dagger:=\bigcup_{\bar r\to0} \bA_K^{(0,\bar r]}$, and}{$\DS \bigcup_{r\to\infty}\TS R\langle\frac{z}{\zeta^r}\rangle$, and}

\comp{$\DS\bB_K^\dagger:=\bigcup_{\bar r\to0} \bB_K^{(0,\bar r]}$ the subfield of $\bB_K$ of \emph{overconvergent elements},}{$\DS \bigcup_{r\to\infty}\TS R\langle\frac{z}{\zeta^r}\rangle[\frac{1}{z}]$.}

\comp{$\bB^\dagger:=\wt\bB^\dagger\cap\bB$.}{\mbox{ }}

\comp{$\phi:\wt\bA^{(0,\bar r]}\isoto\wt\bA^{(0,\bar r/p]}$ is a bicontinuous isomorphism}{$\sigma:C\langle\frac{z}{\zeta^r}\rangle\isoto C\langle\frac{z}{\zeta^{qr}}\rangle$ is a bicontinuous isomorphism}

\comp{$\phi:\bA_K^{(0,\bar r]}\isoto\bA_K^{(0,\bar r/p]}$ is a continuous homomorphism}{$\sigma:R\langle\frac{z}{\zeta^r}\rangle\isoto R\langle\frac{z}{\zeta^{qr}}\rangle$ is a continuous homomorphism}

\comp{Note the inclusions $\bA_K^+\subset\bA_K^{(0,\bar r]}\subset\bA_K$ and the corresponding inclusions for the $\bB$'s or for the ``tilde versions''.}
{Note the inclusions $R\open\dbl z\dbr\subset R\langle\frac{z}{\zeta^r}\rangle\subset R\dbl z\dbr$ and the corresponding inclusions after adjoining $\frac{1}{z}$ or replacing $R$ by $C$.}

\comp{There are natural inclusions\\
$\bA_K^+\subset\bA_K^{(0,\bar r]}\subset\bA_K^\dagger\subset\bA_K$ identifying $\bA_K$ with the $p$-adic completion of $\bA_K^{(0,\bar r]}$ and $\bA_K^\dagger$. The same holds for the $\bB$'s and the ``tilde versions''.}
{There are natural inclusions\\
$R\open\dbl z\dbr\subset R\langle\frac{z}{\zeta^r}\rangle\subset R\dbl z\dbr$ identifying $R\dbl z\dbr$ with the $z$-adic completion of $ R\langle\frac{z}{\zeta^r}\rangle$. The same holds after adjoining $\frac{1}{z}$ or replacing $R$ by $C$.}

\subsection{The Rings $\wt\bB_{\max}$ and $\wt{\bB}_{\rig}$} \label{SectWTBMax}

\comp{Let $\bar r\in\BR_{>0}$ and set $r=1/\bar r$.}{Let $r\in\BR_{>0}$ and set $\bar r=1/r$.}

\comp{Consider the semi-valuation $v^{(0,\bar r]}$ on $\wt\bB^+$.}
{Consider the norm $\|\,.\,\|_r$ on $\CO_C\dbl z\dbr[\frac{1}{z}]$.}

\comp{Let $\wt\bA^+_{\max}$ be the $p$-adic completion of $\wt\bA^+\bigl[\frac{[\tilde p]}{p}\bigr]$ or equivalently its completion with respect to $v^{(0,1]}$. (This ring is usually denoted $\bA_{\max}$.)}
{Let $\CO_C\dbl z,\frac{\zeta}{z}\rangle$ be the $\zeta$-adic completion of $\CO_C\dbl z\dbr[\frac{\zeta}{z}]$ or equivalently the completion with respect to $\|\,.\,\|_1$. Its elements converge on the half open annulus $\{|\zeta|\le|z|<1\}$ and are bounded by $1$.}

\comp{Set $\wt\bB^+_{\max}:=\wt\bA^+_{\max}[\frac{1}{p}]$. It equals the completion of $\wt\bB^+$ with respect to $v^{(0,1]}$ (and is usually denoted $\bB^+_{\max}$.)}{Then $\CO_C\dbl z,\frac{\zeta}{z}\rangle[\frac{1}{z}]$ equals the completion\\
$\DS\bigl\{\,\sum_{i=-\infty}^\infty b_i z^i: b_i\in\CO_C,\lim_{i\to-\infty}|b_i|\,|\zeta|^{i}=0\,\bigr\}$\\
of $\CO_C\dbl z\dbr[\frac{1}{z}]$ with respect to $\|\,.\,\|_1$.\\
Its elements converge on the half open annulus $\{|\zeta|\le|z|<1\}$ and are bounded by $1$ as $|z|\to1$.}

\comp{Set $\wt\bB^+_\rig:=\DS\bigcap_{n\in\BN_0}\phi^n\bigl(\wt\bB^+_{\max}\bigr)$.}
{Set $\CO_C\dbl z,\frac{1}{z}\}:=\DS\bigcap_{n\in\BN_0}\sigma^n\bigl(\TS\CO_C\dbl z,\frac{\zeta}{z}\rangle[\frac{1}{z}]\bigr)=$\\[-1mm]
$\DS\bigl\{\,\sum_{i=-\infty}^\infty b_i z^i: b_i\in\CO_C,\lim_{i\to-\infty}|b_i|\,|\zeta|^{ri}=0\;\,\forall r>0\,\bigr\}$. Its elements converge on the punctured open unit disc $\{0<|z|<1\}$ and are bounded by $1$ as $|z|$ approaches $1$.}

\comp{$\wt\bB^+_\rig$ also equals the Fr\'echet completion of $\wt\bB^+$ with respect to the family of semi-valuations 
$v^{(0,\bar r]}$ for $0<\bar r\le1$.}
{$\CO_C\dbl z,\frac{1}{z}\}$ also equals the Fr\'echet completion of $\CO_C\dbl z\dbr[\frac{1}{z}]$ with respect to the family of norms $\|\,.\,\|_r$ for $1\le r$.}

\comp{The $\phi$-invariants are}{The $\sigma$-invariants are}

\comp{$(\wt\bB^+_{\max})^{\phi=1}=(\wt\bB^+_\rig)^{\phi=1}=\BQ_p$.}
{$\CO_C\dbl z,\frac{\zeta}{z}\rangle[\frac{1}{z}]^{\sigma=1}=\CO_C\dbl z,\frac{1}{z}\}^{\sigma=1}=\BF_q\dpl z\dpr$.}

\comp{The Galois invariants are\\[1mm]
$(\wt\bB^+_{\rig})^{\CG_K}=K_0:=W(k)[\frac{1}{p}]$.}
{In contrast, the Galois invariants are\\[1mm]
$\CO_C\dbl z,\frac{1}{z}\}^{\CG_L}=\CO_{L^p}\dbl z,\frac{1}{z}\}\neq \ell\dpl z\dpr$ making the notion of crystalline Galois representation problematic, see \ref{SectCrystReps}.}

\subsection{The Rings $\bB^+_{K,\max}$ and $\bB^+_{K,\rig}$}

\comp{\mbox{ }}
{Let $R\open$ be an admissible formal $\CO_L$-algebra with $R\open\otimes_{\CO_L}L\cong R$ as in \ref{SectNonPerfect}.}

\comp{Let $\bB^+_{K,\max}$ be the completion of $\bB^+_K$ with respect to the semi-valuation $v^{(0,1]}$.}{Let $R\open\dbl z,\frac{\zeta}{z}\rangle[\frac{1}{z}]$ be the completion of $R\open\dbl z\dbr[\frac{1}{z}]$ with respect to the norm $\|\,.\,\|_1$. It equals}

\comp{\mbox{ }}
{$\DS\bigl\{\,\sum_{i=-\infty}^\infty b_i z^i: b_i\in R\open,\lim_{i\to-\infty}|b_i|\,|\zeta|^{i}=0\,\bigr\}$. Its elements converge on the relative half open annulus $\{|\zeta|\le|z|<1\}$ over $\Spm R$ and are bounded by $1$ as $|z|\to1$.}

\comp{$\bA^+_{K,\max}:=\bigl\{\,x\in\bB^+_{K,\max}: v^{(0,1]}(x)\ge0\,\bigr\}$.}
{$R\open\dbl z,\frac{\zeta}{z}\rangle :=\bigl\{\,x\in R\open\dbl z,\frac{\zeta}{z}\rangle[\frac{1}{z}]:\|x\|_1\le1\,\bigr\}$ equals the $\zeta$-adic completion of $R\open\dbl z\dbr[\frac{\zeta}{z}]$ or equivalently the completion with respect to $\|\,.\,\|_1$.}

\comp{Let $\bB^+_{K,\rig}$ be the Fr\'echet completion of $\bB_K^+$ with respect to the family of semi-valuations 
$v^{(0,\bar r]}$ for $0<\bar r\le1$}
{Let $R\open\dbl z,\frac{1}{z}\}$ be the Fr\'echet completion of $R\open\dbl z\dbr[\frac{1}{z}]$ with respect to the family of norms $\|\,.\,\|_r$ for $1\le r$. It equals\\
$\DS\bigl\{\,\sum_{i=-\infty}^\infty b_i z^i: b_i\in R\open,\lim_{i\to-\infty}|b_i|\,|\zeta|^{ri}=0\;\,\forall r>0\,\bigr\}$. Its elements converge on the relative punctured open unit disc $\{0<|z|<1\}$ over $\Spm R$ and are bounded by $1$ as $|z|\to1$.}

\comp{The $\phi$-invariants are}{If $\Spec R\open$ is connected the $\sigma$-invariants are}

\comp{$(\bB^+_{K,\max})^{\phi=1}=(\bB^+_{K,\rig})^{\phi=1}=\BQ_p$.}
{$R\open\dbl z,\frac{\zeta}{z}\rangle[\frac{1}{z}]^{\sigma=1}=R\open\dbl z,\frac{1}{z}\}^{\sigma=1}=\BF_q\dpl z\dpr$.}

\subsection{The Analogues of $2\pi i$} \label{SectThePeriods}

\comp{The series $t:=\log[\epsilon]:=\DS\sum_{n\in\BN}\TS\frac{(-1)^{n-1}}{n}\pi^n$ converges in $\bB^+_{\BQ_p,\rig}$ and appears as the $p$-adic analogue of $2\pi i$. It satisfies}
{\addtolength{\baselineskip}{1mm}
The product $t_{\SSC -}:=\DS\prod_{n\in\BN_0}\TS\bigl(1-\frac{\zeta^{q^n}}{z}\bigr)$ converges in $\BF_q\dbl\zeta\dbr\dbl z,\frac{1}{z}\}$. Using the notation of \ref{SectCyclo} we let $\BF_q\dbl\zeta\dbr_\infty:=\BF_q\dbl\zeta\dbr[t_n:n\in\BN_0]$ be the valuation ring of $\BF_q\dpl\zeta\dpr_\infty$. Then we define \\[1mm]
$t:=t_{\SSC +}\cdot t_{\SSC -}\in\BF_q\dbl\zeta\dbr_\infty\dbl z,\frac{1}{z}\}$. \es It satisfies }

\comp{$\phi(t)=p\, t$ \es and \es $\gamma(t)=\chi(\gamma)\cdot t$ \\[1mm]
for all $\gamma\in\Gal(\BQ_p^\alg/\BQ_p)$ where $\chi$ is the cyclotomic character from \ref{SectCyclo}.}
{$\sigma(t)=z\, t$ \es and \es $\gamma(t)=\chi(\gamma)\cdot t$ \\[1mm]
for all $\gamma\in\Gal\bigl(\BF_q\dpl\zeta\dpr^\sep/\BFZ\bigr)$ where $\chi$ is the cyclotomic character from \ref{SectCyclo}.}

\comp{Set $\wt\bB_{\max}:=\wt\bB^+_{\max}[\frac{1}{t}]$\;, $\wt\bB_\rig:=\wt\bB^+_\rig[\frac{1}{t}]$\;, \\[1mm]
$\bB_{K,\max}:=\bB^+_{K,\max}[\frac{1}{t}]$\;, and $\bB_{K,\rig}:=\bB^+_{K,\rig}[\frac{1}{t}]$}{For the convergence behavior of the analogues of these rings note that $t$ has simple zeroes precisely at $z=\zeta^{q^i}$ for $i\in\BZ$.}

\subsection{The Robba Ring} \label{SectRobba}

\comp{Let $\bar r,\bar s\in\BR_{>0}$ with $\bar s\le\bar r$ and $r=1/\bar r$, $s=1/\bar s$.}{Let $r,s\in\BR_{>0}$ with $s\ge r$ and $\bar r=1/r$, $\bar s=1/s$.}

\comp{Consider the semi-valuation on $\wt\bB^{(0,r]}$\\
$v^{[\bar s,\bar r]}(x):=\inf\bigl\{\,v^{(0,\bar s]}(x),v^{(0,\bar r]}(x)\,\bigr\}$ and let}
{Consider the norm on $C\langle\frac{z}{\zeta^r}\rangle[\frac{1}{z}]$ or on $R\langle\frac{z}{\zeta^r}\rangle[\frac{1}{z}]$\\
$\|x\|_{[r,s]}:=\sup\bigl\{\,\|x\|_r,\|x\|_s\,\bigr\}$ and let}

\comp{$\wt\bB^{[\bar s,\bar r]}$ (respectively $\bB_K^{[\bar s,\bar r]}$) be the completion of $\wt\bB^{(0,r]}$ (respectively $\bB_K^{(0,r]}$) with respect to $v^{[\bar s,\bar r]}$.}
{$C\langle\frac{z}{\zeta^r},\frac{\zeta^s}{z}\rangle$ (respectively $R\langle\frac{z}{\zeta^r},\frac{\zeta^s}{z}\rangle$) be the completion of $C\langle\frac{z}{\zeta^r}\rangle[\frac{1}{z}]$ (respectively $R\langle\frac{z}{\zeta^r}\rangle[\frac{1}{z}]$) with respect to $\|\,.\,\|_{[r,s]}$. It equals}

\comp{\mbox{ }}
{$\DS\bigl\{\,\sum_{i=-\infty}^\infty b_i z^i: \lim_{i\to\infty}|b_i|\,|\zeta|^{ri} = \lim_{i\to-\infty} |b_i|\,|\zeta|^{si}=0\,\bigr\}$ \\
the ring of Laurent series with coefficients $b_i\in C$ (respectively $b_i\in R$), which converge on the (relative) annulus $\{|\zeta|^s\le|z|\le|\zeta|^r\}$.}

\comp{Let $\wt\bB^{]0,\bar r]}$ (respectively $\bB_K^{]0,\bar r]}$) be the Fr\'echet completion of $\wt\bB^{(0,r]}$ (respectively $\bB_K^{(0,r]}$) with respect to the family of semi-valuations $v^{[\bar s,\bar r]}$ for all $0<\bar s\le\bar r$. Then}
{Let $C\langle\frac{z}{\zeta^r},\frac{1}{z}\}$ (respectively $R\ancon[r]$) be the Fr\'echet completion of $C\langle\frac{z}{\zeta^r}\rangle[\frac{1}{z}]$ (respectively $R\con[r]$) with respect to the family of norms $\|\,.\,\|_{[r,s]}$ for all $s\ge r$. Then}

\comp{$\DS\wt\bB^{]0,\bar r]}=\bigcap_{\bar s\to0}\wt\bB^{[\bar s,\bar r]}$ \es and \es $\DS\bB_K^{]0,\bar r]}=\bigcap_{\bar s\to0}\bB_K^{[\bar s,\bar r]}$.}
{$C\langle\frac{z}{\zeta^r},\frac{1}{z}\}=\DS\bigcap_{s\to\infty}\TS C\langle\frac{z}{\zeta^r},\frac{\zeta^s}{z}\rangle=$}

\comp{The functions $f_i$ from \ref{SectWittPower} extend by continuity to $\wt\bB^{[\bar s,\bar r]}$ and $\wt\bB^{]0,\bar r]}$, and the sum $\sum_{i=-\infty}^\infty p^i[f_i(x)]$ converges to $x$ in $\wt\bB^{[\bar s,\bar r]}$, respectively in $\wt\bB^{]0,\bar r]}$. This is \emph{not true} for $\bB_K^{[\bar s,\bar r]}$ and $\bB_K^{]0,\bar r]}$ since $\bE_K$ is not perfect and so there is no Teich\-m\"uller map $\bE_K\to\bA_K$.}
{$\DS\bigl\{\,\sum_{i=-\infty}^\infty b_i z^i:\lim_{i\to\pm\infty}\! |b_i|\,|\zeta|^{si}=0\text{ for all }s\ge r\,\bigr\}$ \\
equals the ring of Laurent series with coefficients $b_i\in C$, which converge on the punctured disc $\{0<|z|\le|\zeta|^r\}$. This is \emph{also true} if $C$ is replaced by $R$.}

\comp{The rings $\wt\bB^{[\bar s,\bar r]}$ and $\bB_K^{[\bar s,\bar r]}$ are principal ideal domains, the rings  $\wt\bB^{]0,\bar r]}$ and $\bB_K^{]0,\bar r]}$ are Bezout domains (that is, every finitely generated ideal is principal).}
{The rings $C\langle\frac{z}{\zeta^r},\frac{\zeta^s}{z}\rangle$ and $L\langle\frac{z}{\zeta^r},\frac{\zeta^s}{z}\rangle$ are principal ideal domains, the rings $C\ancon[r]$ and $L\ancon[r]$ are Bezout domains.}

\comp{For $[\bar s',\bar r']\subset[\bar s,\bar r]$ one has inclusions with dense image \raisebox{0em}[1.2em]{$\wt\bB^{(0,\bar r]}\subset\wt\bB^{]0,\bar r]}\subset\wt\bB^{[\bar s,\bar r]}\subset\wt\bB^{[\bar s',\bar r']}$} and the same for the $\bB_K$-versions.}
{For $[r',s']\subset[r,s]$ one has the dense inclusions \\ 
$C\langle\frac{z}{\zeta^r}\rangle[\frac{1}{z}]\subset C\langle\frac{z}{\zeta^r},\frac{1}{z}\}\subset C\langle\frac{z}{\zeta^r},\frac{\zeta^s}{z}\rangle\subset C\langle\frac{z}{\zeta^{r'}},\frac{\zeta^{s'}}{z}\rangle$ and the same if $C$ is replaced by $L$ or $R$.}

\comp{Set $\DS\wt\bB^\dagger_\rig:=\bigcup_{\bar r\to 0}\wt\bB^{]0,\bar r]}$ and\\[1mm]
$\DS\bB^\dagger_{K,\rig}:=\bigcup_{\bar r\to 0}\bB_K^{]0,\bar r]}$. The later ring is called the \emph{Robba ring} associated with $K$.}
{$\DS\bigcup_{r\to\infty}\TS C\langle\frac{z}{\zeta^r},\frac{1}{z}\}$ (respectively $\DS\bigcup_{r\to\infty}\TS R\langle\frac{z}{\zeta^r},\frac{1}{z}\}$) is the ring of Laurent series which converge on some (relative) punctured disc $\{0<|z|\le|\zeta|^r\}$ with small enough radius $|\zeta|^r$.}

\comp{$\phi$ induces isomorphisms of topological rings\\[2mm]
$\phi:\wt\bB^{[\bar s,\bar r]}\isoto\wt\bB^{[\bar s/p\,,\,\bar r/p]}$,}
{$\sigma$ induces isomorphisms of topological rings \\[2mm]
$\sigma:C\langle\frac{z}{\zeta^r},\frac{\zeta^s}{z}\rangle\isoto C\langle\frac{z}{\zeta^{qr}},\frac{\zeta^{qs}}{z}\rangle$,}

\comp{$\phi:\wt\bB^{]0,\bar r]}\isoto\wt\bB^{]0\,,\,\bar r/p]}$,}
{$\sigma:C\langle\frac{z}{\zeta^r},\frac{1}{z}\}\isoto C\langle\frac{z}{\zeta^{qr}},\frac{1}{z}\}$,}

\comp{$\phi:\wt\bB^\dagger_\rig\isoto\wt\bB^\dagger_\rig$,}
{$\sigma:\DS\bigcup_{r\to\infty}\TS C\langle\frac{z}{\zeta^r},\frac{1}{z}\}\isoto \DS\bigcup_{r\to\infty}\TS C\langle\frac{z}{\zeta^r},\frac{1}{z}\}$,}

\comp{and homomorphisms of topological rings\\[2mm]
$\phi:\bB_K^{[\bar s,\bar r]}\longto\bB_K^{[\bar s/p\,,\,\bar r/p]}$,}
{and homomorphisms of topological rings \\[2mm]
$\sigma:R\langle\frac{z}{\zeta^r},\frac{\zeta^s}{z}\rangle\longto R\langle\frac{z}{\zeta^{qr}},\frac{\zeta^{qs}}{z}\rangle$,}

\comp{$\phi:\bB_K^{]0,\bar r]}\longto\bB_K^{]0\,,\,\bar r/p]}$,}
{$\sigma:R\langle\frac{z}{\zeta^r},\frac{1}{z}\}\longto R\langle\frac{z}{\zeta^{qr}},\frac{1}{z}\}$,}

\comp{$\phi:\bB^\dagger_{K,\rig}\longto\bB^\dagger_{K,\rig}$.}
{$\sigma:\DS\bigcup_{r\to\infty}\TS R\langle\frac{z}{\zeta^r},\frac{1}{z}\}\longto \DS\bigcup_{r\to\infty}\TS R\langle\frac{z}{\zeta^r},\frac{1}{z}\}$.}

\comp{There are exact sequences of rings
\[
\begin{array}{@{0\;\to\es}l@{\es\to\es}l@{\,\oplus\,}l@{\es\to\es}l@{\es\to\;0}l}
\wt\bB^+ & \wt\bB^+_\rig & \wt\bB^{(0,\bar r]} & \wt\bB^{]0,\bar r]} &, \\[2mm]
\wt\bB^+ & \wt\bB^+_{\max} & \wt\bB^{(0,1]} & \wt\bB^{[1,1]} &,
\end{array}
\]
and the same for the $\bB_K$-versions.}
{There are exact sequences of rings
\[
\begin{array}{l@{\;\hookrightarrow\;}l}
\CO_C\dbl z\dbr[\frac{1}{z}] & \CO_C\dbl z,\frac{1}{z}\} \oplus C\langle\frac{z}{\zeta^r}\rangle[\frac{1}{z}] \shortonto C\langle\frac{z}{\zeta^r},\frac{1}{z}\}, \\[2.6mm]
\CO_C\dbl z\dbr[\frac{1}{z}] & \CO_C\dbl z,\frac{\zeta}{z}\rangle[\frac{1}{z}] \oplus C\langle\frac{z}{\zeta}\rangle[\frac{1}{z}] \shortonto C\langle\frac{z}{\zeta},\frac{\zeta}{z}\rangle,
\end{array}
\]
and the same if $\CO_C$ and $C$ are replaced by $R\open$ and $R$.}

\subsection{The Field $\bB_\dR$} \label{SectBdR}

\comp{The surjective homomorphism \\[1mm]
$\theta:\wt\bB^+\to C\,,\;x\mapsto\DS\sum_{i\in\BZ}p^i f_i(x)^{(0)}$ extends by continuity to $\theta:\wt\bB^{[1,1]}\to C$. It has}
{Let $\theta:C\langle\frac{z}{\zeta},\frac{\zeta}{z}\rangle\to C\,,\;\DS\sum_{i=-\infty}^\infty b_iz^i \mapsto \sum_{i=-\infty}^\infty b_i\zeta^i$.\\
 It is a surjective homomorphism with}

\comp{$\ker\theta=(p-[\tilde p])\wt\bB^{[1,1]}=t\,\wt\bB^{[1,1]}$.\\[1mm]
However, note that the ideal $t\,\wt\bB^+_\rig$ of $\wt\bB^+_\rig$ is not maximal as opposed to $t\,\wt\bB^{[1,1]}\subset\wt\bB^{[1,1]}$.}
{$\ker\theta=(z-\zeta)C\langle\frac{z}{\zeta},\frac{\zeta}{z}\rangle=t\,C\langle\frac{z}{\zeta},\frac{\zeta}{z}\rangle$.\\[1mm]
However, note that $t$ also has other zeroes outside $\{|z|=|\zeta|\}$.}

\comp{Set $\bB^+_\dR:=\invlim[n]\wt\bB^{[1,1]}/(\ker\theta)^n$. It is a complete discrete valuation ring with residue field $C$ and uniformizing parameter $p-[\tilde p]$ or $t$. There is \emph{no} continuous section $C\to\bB^+_\dR$.}
{Let $C\dbl z-\zeta\dbr=\invlim[n] C\langle\frac{z}{\zeta},\frac{\zeta}{z}\rangle/(\ker\theta)^n$ be the power series ring over $C$ in the ``variable'' $z-\zeta$. It is a complete discrete valuation ring with residue field $C$ and uniformizing parameter $z-\zeta$ or $t$ \emph{with }canonical section $C\to C\dbl z-\zeta\dbr$. It is the complete local ring at the point $\{z=\zeta\}=\ker\theta$.}

\comp{$\bB_\dR:=\bB^+_\dR[\frac{1}{t}]$ is Fontaine's \emph{$p$-adic period field}.}
{$C\dpl z-\zeta\dpr=C\dbl z-\zeta\dbr[\frac{1}{t}]$ is its analogue.}

\comp{It is filtered by putting $Fil^n\bB_\dR:=t^n\bB^+_\dR$.}
{It is filtered by putting \\
$Fil^n C\dpl z-\zeta\dpr:=t^n C\dbl z-\zeta\dbr$ for $n\in\BZ$.}

\comp{$\bB_\dR$ is a successive extension of the Galois modules $C(n)$ for $n\in\BZ$ as follows (see \ref{SectTateTwists})
\[
0\to t^{n+1}\bB^+_\dR\to t^n\bB^+_\dR\to C(n)\to 0\,.
\]
}
{$C\dpl z-\zeta\dpr$ is a successive extension of the Galois modules $C(n)\cong C$ for $n\in\BZ$ as follows (see \ref{SectTateTwists})
\[
0\to t^{n+1}C\dbl z-\zeta\dbr\to t^nC\dbl z-\zeta\dbr\to C\to 0\,.
\]
}

\comp{This implies $(\bB_\dR)^{\CG_K}=K$.}
{This implies the unfavorable fact that \\
$C\dpl z-\zeta\dpr^{\CG_L}=L^p\dpl z-\zeta\dpr\supsetneq L$ which spoils the use of $G_L$-representations; see Section~\ref{ChaptGaloisRep}.}

\comp{$\wt\bB^+_\dR$ naturally contains $\wt\bB^{[1,1]}$ and if $p^{-n}\in[\bar s,\bar r]$ for some $n\in\BZ$ then this induces an inclusion $\phi^{-n}:\wt\bB^{[\bar s,\bar r]}\hookrightarrow\bB^+_\dR$.}{$C\dbl z-\zeta\dbr$ naturally contains $C\langle\frac{z}{\zeta},\frac{\zeta}{z}\rangle$ and if $q^{n}\in[r,s]$ for some $n\in\BZ$ this induces an inclusion $\sigma^{-n}:C\langle\frac{z}{\zeta^r},\frac{\zeta^s}{z}\rangle\hookrightarrow C\dbl z-\zeta\dbr$.}

\comp{By \cite{Colmez02} the following sequence is exact} 
{It corresponds to the exact sequence}

\comp{$0\to\BQ_p\to(\wt\bB_{\max}){}^{\phi=1}\to\bB_\dR/\bB^+_\dR\to0$.}
{$0\to\BF_q\dpl z\dpr\to\CO_C\dbl z,\frac{\zeta}{z}\rangle[\frac{1}{z}]^{\sigma=1}\to\DS\frac{C\dpl z-\zeta\dpr}{C\dbl z-\zeta\dbr}\to0$.}

\comp{Note that one can construct relative versions of $\bB^+_\dR$ and $\bB_\dR$ over affinoid $K$-algebras; see \cite{Brinon,Wintenberger94}}{There are obvious $R$-rational versions $R\dbl z-\zeta\dbr$ and $R\dbl z-\zeta\dbr[\frac{1}{z-\zeta}]$ of these rings over an arbitrary affinoid $L$-algebra $R$.}

\section{$p$-Divisible Groups, Local Shtuka and Crystals}\label{ChaptLocalShtuka}

\subsection{Barsotti-Tate Groups} \label{SectBTGroups}

\comp{Let $S$ be a scheme and $h\ge0$ be an integer. A \emph{Barsotti-Tate group} or \emph{$p$-divisible group of height $h$} over $S$ is an inductive system 
\[
G \es = \es (G_1\xrightarrow{i_1}G_2\xrightarrow{i_2}G_3\xrightarrow{i_3}\ldots)
\]
where for each $n\ge1$
\begin{itemize}
\item 
$G_n$ is a finite commutative group scheme over $S$ of order $p^{nh}$,
\item 
the sequence of group schemes over $S$ is exact
\[
0\to G_n\xrightarrow{i_n}G_{n+1}\xrightarrow{p^n}G_{n+1}\,.
\]
\end{itemize}
Obviously multiplication with $p$ on $G_n$ induces on the Lie algebra of $G_n$ multiplication by the scalar $p$, compare \ref{SectTwoRolesOfP}.}
{Let $S$ be an $\BF_q[z]$-scheme and denote the image of $z$ in $\CO_S$ by $\zeta$.
Let $d,h\ge0$ be integers. A \emph{divisible local Anderson module of height $h$ and dimension $d$} over $S$ is an inductive system of finite $\BF_q\dbl z\dbr$-module schemes over $S$
\[
G \es = \es (G_1\xrightarrow{i_1}G_2\xrightarrow{i_2}G_3\xrightarrow{i_3}\ldots)
\]
where for each $n\ge 1$
\begin{itemize}
\itemsep0ex
\item \label{DefZDivGpAxiom1}
the $\BF_q$-module scheme $G_n$ can be embedded into an $\BF_q$-vector group scheme over $S$,
\item \label{DefZDivGpAxiom2}
the order of $G_n$ is $q^{hn}$,
\item \label{DefZDivGpAxiom3}
the following sequence of $\BF_q\dbl z\dbr$-module schemes over $S$ is exact
\[
0\to G_n\xrightarrow{i_n}G_{n+1}\xrightarrow{z^n}G_{n+1}\,,
\]
\item \label{DefZDivGpAxiom4}
$(z-\zeta)^d=0$ on $\Lie G_n$\,, 
\item \label{DefZDivGpAxiom5}
$\DS d \es = \es\max_{ s\in S,\, n\ge1}\{\,\dim_{\kappa(s)}\bigl(\Lie G_n\otimes_{\CO_S}\kappa(s)\bigr)\,\}$\,.
\end{itemize}
So $z$ does not need to act on $\Lie G_n$ as the scalar $\zeta$, compare \ref{SectTwoRolesOfP}}

\subsection{Local Shtuka} \label{SectLocalShtuka}

\comp{Let $k\supset\BF_p$ be a perfect field. A \emph{Dieudonn\'e crystal} over $k$ is a pair $(M,F_M)$ where $M$ is a finite free $W(k)$-module and $F_M: M\to M$ is a $\phi$-linear endomorphism with $pM\subset F_M(M)$.\\[2mm]
Equivalently we can set $\phi^\ast M:=M\otimes_{W(k),\phi}W(k)$ and linearize $F_M$ to a homomorphism of $W(k)$-modules $F_M:\phi^\ast M\to M$ which satisfies $p=0$ on $\coker F_M$.}
{Let $S$ be a formal scheme over $\BF_q\dbl\zeta\dbr$. A \emph{local shtuka} of rank $n$ over $S$ is a pair $(M,F_M)$ consisting of a sheaf of $\CO_S\dbl z\dbr$-modules on $S$ and an isomorphism $F_M:\sigma^\ast M[\frac{1}{z-\zeta}] \isoto M[\frac{1}{z-\zeta}]$, where $\sigma^\ast M:=M\otimes_{\CO_S\dbl z\dbr,\sigma}\CO_S\dbl z\dbr$, such that the following conditions hold:
\begin{itemize}
\item
locally for the Zariski topology on $S$, $M$ is a free $\CO_S\dbl z\dbr$-module of rank $n$,
\item
there exists an integer $e$ such that $F_M(\sigma^\ast M) \subset (z-\zeta)^{-e} M$ and the quotient $(z-\zeta)^{-e} M/F_M(\sigma^\ast M)$ is locally free and coherent as an $\CO_S$-module.
\end{itemize}
A local shtuka is called \emph{effective} if $F_M$ is actually a morphism $\sigma^\ast M\to M$.}

\comp{Then the functor which assigns to a Barsotti-Tate group its contravariant Dieudonn\'e module is an anti-equivalence between the category of Barsotti-Tate groups over $k$ and the category of Dieudonn\'e crystals over $k$, see \cite[p.\ 71]{Demazure}.}
{Then the category of divisible local Anderson modules over $S$ is anti-equivalent to the category of effective local shtuka over $S$; see \cite{Crystals}.}

\comp{For a Barsotti-Tate group over $k$ the properties of being \'etale or connected and the description of isogenies reflect in its Dieudonn\'e module.}
{For a divisible local Anderson module the properties of being \'etale or having connected fibers, and the description of isogenies can be read off from its associated local shtuka.}

\comp{Barsotti-Tate groups arise from abelian varieties as their subgroups of $p$-power torsion. They are of most interest when $p$ equals the characteristic of the ground field.}
{Also divisible local Anderson modules arise from global objects like Drinfeld-modules \cite{Drinfeld} and abelian $t$-modules \cite{Anderson} as their subgroups of $z$-power torsion. Similarly local shtuka arise from global objects like shtuka \cite{Laumon97}, $t$-motives \cite{Anderson}, or abelian sheaves \cite{AbSh} by completing with respect to $z$.}

\commt{The parallel between Barsotti-Tate groups or more generally $F$-crystals and local shtuka is close. It ranges from the classification over algebraically closed fields (see \ref{SectDM1}), the behavior of Newton and Hodge polygons, like the Grothendieck-Katz Specialization Theorem, to their deformation theory; see Katz~\cite{Katz}, Grothendieck~\cite{Grothendieck}, Messing~\cite{Messing}, Hartl~\cite{Crystals} and \cite[\S\S 6--8]{AbSh}.}

\subsection{Tate Modules} \label{SectTateModules}

\comp{Let $K$ be a complete discretely valued field extension of $\BQ_p$ with perfect residue field $k$.}
{Let $L$ be a field extension of $\BFZ$ which is complete with respect to an absolute value $|\,.\,|:L\to\BR_{\geq0}$ extending the absolute value on $\BFZ$.}

\comp{The \emph{Tate module} of a Barsotti-Tate group $G=\dirlim(G_n,i_n)$ over $\CO_K$ is the $\BZ_p[\CG_K]$-module 
\[
T_pG\;:=\;\invlim[n]\bigl(G_n(K^\alg),p\bigr)\,.
\]}
{The \emph{Tate module} of a divisible local Anderson module $G=\dirlim(G_n,i_n)$ over $\CO_L$ is the $\BF_q\dbl z\dbr[\CG_L]$-module
\[
T_zG\;:=\;\invlim[n]\bigl(G_n(L^\sep),z\bigr)\,.
\]}

\comp{\mbox{ }}
{The \emph{Tate module} $T_z(M,F_M)$ of a local shtuka $(M,F_M)$ over $\CO_L$ is the $\BF_q\dbl z\dbr$-dual of $(M,F_M)^{F=1}(L^\sep)$ which by definition equals
\[
\bigl\{\,m\in M\otimes_{\CO_L\dbl z\dbr}L^\sep\dbl z\dbr : F_M(\sigma^\ast m)=m\,\bigr\}\,.
\]}

\comp{There is a description of $T_pG$ in terms of the Dieudonn\'e crystal associated with $G$, see Fontaine~\cite[\S V.1]{Fontaine77}.}
{If $M(G)$ is the local shtuka over $\BF_q\dbl\zeta\dbr$ associated with $G$ then the $\BF_q\dbl z\dbr[\CG_L]$-modules $T_zG$ and $T_zM(G)$ are canonically isomorphic; see \cite{Crystals}.}

\subsection{Lubin-Tate Formal Groups} \label{SectLubinTate}

\comp{Let $G$ be the Lubin-Tate formal group over $\BZ_p$ on which $p$ acts as $x\mapsto (1+x)^p-1$, see \cite{LT}.}
{Let $G$ be the formal additive group over $\BF_q\dbl\zeta\dbr$ on which we let $z$ act by $x\mapsto \zeta x+x^q$. Then $G$ is a }

\comp{Then $G$ is a Barsotti-Tate group of height $1$ over $\BZ_p$. It is isomorphic over $\BZ_p$ to the Lubin-Tate formal group on which $p$ acts as $x\mapsto px+x^p$.}
{divisible local Anderson module of height $1$ and dimension $1$. If we identify $z$ with $\zeta$ then $G$ is the Lubin-Tate formal group over $\BF_q\dbl \zeta\dbr$ on which $\zeta$ acts as $x\mapsto\zeta x+x^q$. The local shtuka associated with $G$ is $M(G)=\bigl(\BF_q\dbl\zeta\dbr\dbl z\dbr,(z-\zeta)\cdot\sigma\bigr)$. }

\comp{The Tate module of $G$ is generated over $\BZ_p$ by 
\[
(\epsilon^{(n)})_{n\in\BN}\;\in\;T_pG\;=\;\invlim[n]\bigl(G_n(K^\alg),p\bigr)\,,
\]
where $\epsilon^{(n)}$ was defined in \ref{SectCyclo}. Thus $\CG_{\BQ_p}$ acts on $T_pG$ through the cyclotomic character from \ref{SectCyclo}.}
{For the Tate module we obtain 
\[
T_zG\;=\;T_zM(G)\;=\;\BF_q\dbl z\dbr\cdot t_{\SSC +}\,,
\]
where $t_{\SSC +}$, as defined in \ref{SectCyclo}, is viewed as the map 
\[
M(G)^{F=1}\bigl(\BFZ^\sep\bigr)\to\BF_q\dbl z\dbr\,,\quad (t_{\SSC +})^{-1}\mapsto 1\,.
\]
Thus $\CG_\BFZ$ acts on $T_zG$ through the cyclotomic character from \ref{SectCyclo}.}

\subsection{Crystals and Isocrystals} \label{SectCrystals}

\comp{Let $k$ be a perfect field containing $\BF_p$. An \emph{$F$-(iso)crystal} over $k$ is a pair $(D,F_D)$ where
\begin{itemize}
\item 
$D$ is a finite free module over $W(k)$ (respectively over $W(k)[\frac{1}{p}]$\;),
\item 
$F_D:\phi^\ast D\to D$ is a homomorphism of $W(k)$-modules with $p$-torsion cokernel (respectively which is an isomorphism).
\end{itemize}
Every Dieudonn\'e crystal over $k$ is an $F$-crystal.}
{Let $\ell\supset\BF_q$ be a field which is an $\BF_q\dbl z\dbr$-algebra in which the image $\zeta$ of $z$ is zero. A \emph{$z$-(iso)crystal} over $\ell$ is a pair $(D,F_D)$ where
\begin{itemize}
\item 
$D$ is a finite free module over $\ell\dbl z\dbr$ (respectively over $\ell\dpl z\dpr$\;),
\item 
put $\sigma^\ast D:=D\otimes_{\ell\dbl z\dbr,\sigma}\ell\dbl z\dbr$, then $F_D:\sigma^\ast D\to D$ is a homomorphism of $\ell\dbl z\dbr$-modules with $z$-torsion cokernel (respectively which is an isomorphism).
\end{itemize}
So a $z$-crystal is nothing but an effective local shtuka over $\ell$.}

\subsection{The Dieudonn\'e-Manin Classification} \label{SectDM1}

\comp{For integers $r,d$ with $r>0$ and $(r,d)=1$ consider the $F$-isocrystal $D_{d,r}=\bigl(W(k)[{\TS\frac{1}{p}}]^{\oplus r},F_D\bigr)$ over $k$ with
\[
F_D=\left(
\raisebox{3.6ex}{$
\xymatrix @C=0.3pc @R=0.3pc {
0\ar@{.}[ddrr] & & p^{-d} \\
1\ar@{.}[dr]\\
& 1 & 0
}$}\right)\!\cdot\phi\es.
\]}
{For integers $r,d$ with $r>0$ and $(r,d)=1$ consider the $z$-isocrystal $D_{d,r}=\bigl(\ell\dpl z\dpr^{\oplus r},F_D\bigr)$ over $\ell$ with
\[
F_D=\left(
\raisebox{3.6ex}{$
\xymatrix @C=0.3pc @R=0.3pc {
0\ar@{.}[ddrr] & & z^{-d} \\
1\ar@{.}[dr]\\
& 1 & 0
}$}\right)\!\cdot\sigma\es.
\]}

\comp{If $k$ is algebraically closed every $F$-isocrystal over $k$ is isomorphic to a direct sum $\bigoplus_i D_{d_i,r_i}$ for uniquely determined pairs of integers $(r_i,d_i)$ up to permutation; see Manin~\cite{Manin}.}
{If $\ell$ is algebraically closed every $z$-isocrystal over $\ell$ is isomorphic to a direct sum $\bigoplus_i D_{d_i,r_i}$ for uniquely determined pairs of integers $(r_i,d_i)$ up to permutation; see Laumon~\cite[\S B.1]{Laumon}.}

\subsection{Filtered Isocrystals} \label{SectFilteredIsocrystals}

\comp{Let $K$ be a complete discretely valued extension of $\BQ_p$ with perfect residue field $k$ and set $K_0:=W(k)[\frac{1}{p}]$. }
{Let $L$ be as in \ref{SectBTGroups}. Let $\CO_L$ be its valuation ring and let $\ell$ be its residue field. Assume that there is a fixed section $\ell\hookrightarrow\CO_L$ of the residue map $\CO_L\to\ell$. This induces in particular a homomorphism $\ell\dpl z\dpr\to L\dbl z-\zeta\dbr$ into the power series ring over $L$ sending $z$ to $z=\zeta+(z-\zeta)$. }

\comp{A \emph{filtered isocrystal} $\ul D=(D,F_D,Fil^\bullet D_K)$ over $K$ consists of
\begin{itemize}
\item 
$(D,F_D)$ an $F$-isocrystal over $k$ and
\item 
$Fil^\bullet D_K$ a decreasing separated exhaustive filtration of $D_K:=D\otimes_{K_0}K$ by $K$-subspaces. 
\end{itemize}
One does not require any compatibility between $F_D$ and $Fil^\bullet$.}
{A \emph{filtered isocrystal} $\ul D=(D,F_D,\Fq_D)$ over $L$ consists of
\begin{itemize}
\item 
$(D,F_D)$ a $z$-isocrystal over $\ell$ and
\item 
$\Fq_D$ an $L\dbl z-\zeta\dbr$-lattice inside the $L\dpl z-\zeta\dpr$-vector space $\sigma^\ast D\otimes_{\ell\dpl z\dpr}L\dpl z-\zeta\dpr$.
\end{itemize}
$\Fq_D$ is called a \emph{Hodge-Pink structure} on $(D,F_D)$. One also sets $\Fp_D:=\sigma^\ast D\otimes_{\ell\dpl z\dpr}L\dbl z-\zeta\dbr$. For the significance of $L\dbl z-\zeta\dbr$ and its analogue in mixed characteristic see \ref{SectBdR}.}

\comp{The filtration $Fil^\bullet D_K$ defines a $G_K$-stable $\bB^+_\dR$-lattice $\Fq_D$ (see \ref{SectBdR}) in $\Fp_D\otimes_{\bB_\dR^+}\bB_\dR$ where $\Fp_D:=D_K\otimes_K \bB_\dR$, by setting
\[
\Fq_D\es:=\es Fil^0(D_K\otimes_K \bB_\dR)\,.
\]
Conversely any such lattice determines a filtration by
\[
Fil^i D_K\es=\es\Bigl(\bigl(\Fp_D\cap t^i\Fq_D\bigr)/\bigl(t\,\Fp_D\cap t^i\Fq_D\bigr)\Bigr)^{G_K}\,.
\]
This defines a 1-1-correspondence between filtrations and $G_K$-stable lattices. }{Every Hodge-Pink structure determines in particular a \emph{Hodge-Pink filtration} $Fil^\bullet$ on the $L$-vector space $D_L:=\Fp_D/(z-\zeta)\Fp_D$ by letting $Fil^i D_L$ be the subspace
\[
\Bigl(\Fp_D\;\cap\;(z-\zeta)^i\Fq_D\Bigr)\,/\,\Bigl((z-\zeta)\Fp_D\;\cap\;(z-\zeta)^i\Fq_D\Bigr)
\]
of $D_L$. Note that $z$ acts on $D_L$ as the scalar $\zeta$. However, as was observed by Pink~\cite{Pink}, the fact that $z$ and $\zeta$ both play part of the role of $p$ makes it necessary to consider Hodge-Pink structures instead of only Hodge-Pink filtrations to get a reasonable category. See \cite[Remark 2.2.3]{HartlPSp} for a detailed discussion of this phenomenon.}

\comp{The \emph{Hodge-Tate weights} of $\ul D$ are the integers $h$ for which $Fil^{-h}D_K\ne Fil^{-h+1}D_K$ or equivalently, the elementary divisors of the $\bB^+_\dR$-lattice $\Fq_D$ relative to $\Fp_D$. Note that the definition differs by a minus sign from the corresponding definition in equal characteristic.}
{The \emph{Hodge-Pink weights} of $(D,F_D,\Fq_D)$ are the elementary divisors of $\Fp_D$ relative to $\Fq_D$. More precisely if $e$ is a large enough integer such that $(z-\zeta)^e\Fp_D\subset \Fq_D$ and
\[
\Fq_D/(z-\zeta)^e\Fp_D\es\cong\es\bigoplus_{i=1}^n L\dbl z-\zeta\dbr/(z-\zeta)^{w_i+e}
\]
then the Hodge-Pink weights are the integers $w_1,\ldots,w_n$.}

\subsection{Weak Admissibility} \label{SectWA}

\comp{Let $\ul D=(D,F_D,Fil^\bullet D_K)$ be a filtered isocrystal of rank $n$ over $K$ and define}
{Let $\ul D=(D,F_D,\Fq_D)$ be a filtered isocrystal of rank $n$ over $L$ and define}

\comp{$t_N(\ul D):=\ord_p(\det F_D)$ the \emph{Newton slope} and}
{$t_N(\ul D):=\ord_z(\det F_D)$ the \emph{Newton slope} and}

\comp{$t_H(\ul D):=\DS\sum_{i\in\BZ}i\cdot\dim_K gr^i_{Fil^\bullet} D_K=$\\[1mm]
\mbox{ }\qquad\quad$=\max\{\,i\in\BZ:Fil^i(\wedge^n D_K)=\wedge^n D_K\,\}$ the \emph{Hodge slope} of $\ul D$, which also equals the negative of the sum of the Hodge-Tate weights counted with multiplicity.}
{$t_H(\ul D):=\DS\sum_{i\in\BZ}i\cdot\dim_L gr^i_{Fil^\bullet} D_L$ the \emph{Hodge slope} of $\ul D$, which also equals the integer $e$ such that $\wedge^n\Fq_D=(z-\zeta)^{-e}\wedge^n\Fp_D$ or the sum of the Hodge-Pink weights counted with multiplicity.}

\comp{$\ul D$ is called \emph{weakly admissible} if $t_H(\ul D)=t_N(\ul D)$ and $t_H(\ul D')\leq t_N(\ul D')$ for any subobject}
{$\ul D$ is called \emph{weakly admissible} if $t_H(\ul D)=t_N(\ul D)$ and $t_H(\ul D')\leq t_N(\ul D')$ for any subobject}

\comp{$\ul D'=(D',F_D|_{D'},Fil^\bullet D'_K)\subset\ul D$, where $D'\subset D$ is an $F_D$-stable $K_0$-subspace and $Fil^\bullet D'_K$ is the induced filtration on $D'_K=D'\otimes_{K_0}K$.}
{$\ul D'=(D',F_D|_{D'},\Fq_{D'})\subset\ul D$, where $D'\subset D$ is an $F_D$-stable $\ell\dpl z\dpr$-subspace and the lattice $\Fq_{D'}$ equals $\Fq_D\cap\sigma^\ast D'\otimes_{\ell\dpl z\dpr}L\dpl z-\zeta\dpr$.}

\section{The Slope Filtration Theorem for Frobenius Modules}

\subsection{Frobenius Modules}

\comp{Let $K$ be as in Section~\ref{ChaptFontainesRings} and set $\CR=\bB^\dagger_{K,\rig}$, or $\CR=\wt\bB^\dagger_\rig$, or $\CR=\bB^\dagger_K$, see \ref{SectOverconRings} and \ref{SectRobba}.}
{Let $L$, $C$, and $R$ be as in Section~\ref{ChaptFontainesRings} and set $\CR=\DS\bigcup_{r\to\infty}R\ancon[r]$ or $\CR=\DS\bigcup_{r\to\infty}R\con[r]$, see \ref{SectOverconRings} and \ref{SectRobba}. Note that we allow $R=L=C$.}

\comp{A \emph{$\phi$-module} over $\CR$ consists of a finite free $\CR$-module $M$ and an isomorphism $F_M:\phi^\ast M\isoto M$ where $\phi^\ast M:=M\otimes_{\CR,\phi}\CR$.}
{A \emph{$\sigma$-module} over $\CR$ consists of a finite free $\CR$-module $M$ and an isomorphism $F_M:\sigma^\ast M\isoto M$ where $\sigma^\ast M:=M\otimes_{\CR,\sigma}\CR$. }

\subsection{Dieudonn\'e-Manin Decompositions} 

\comp{For integers $r,d$ with $r>0$ and $(r,d)=1$ consider the $\phi$-module $M_{d,r}=\bigl(\CR^{\oplus r},F_M\bigr)$ over $\CR$ with
\[
F_M=\left(
\raisebox{3.6ex}{$
\xymatrix @C=0.3pc @R=0.3pc {
0\ar@{.}[ddrr] & & p^{-d} \\
1\ar@{.}[dr]\\
& 1 & 0
}$}\right)\!\cdot\phi\es.
\]}
{For integers $r,d$ with $r>0$ and $(r,d)=1$ consider the $\sigma$-module $M_{d,r}=\bigl(\CR^{\oplus r},F_M\bigr)$ over $\CR$ with
\[
F_M=\left(
\raisebox{3.6ex}{$
\xymatrix @C=0.3pc @R=0.3pc {
0\ar@{.}[ddrr] & & z^{-d} \\
1\ar@{.}[dr]\\
& 1 & 0
}$}\right)\!\cdot\sigma\es.
\]}

\comp{Every $\phi$-module over $\wt\bB^\dagger_\rig$ is isomorphic to a direct sum $\bigoplus_i M_{d_i,r_i}$ for uniquely determined pairs of integers $(r_i,d_i)$ up to permutation; see Kedlaya~\cite[Theorem 4.5.7]{Kedlaya}.}
{Every $\sigma$-module over $\bigcup_r C\ancon[r]$ is isomorphic to a direct sum $\bigoplus_i M_{d_i,r_i}$ for uniquely determined pairs of integers $(r_i,d_i)$ up to permutation; see Hartl-Pink~\cite[Theorem 11.1]{HP}.}

\comp{A $\phi$-module $M$ over $\bB^\dagger_{K,\rig}$ is called \emph{(isoclinic of) slope $\lambda$} if $M\otimes_{\bB^\dagger_{K,\rig}}\wt\bB^\dagger_\rig$ is isomorphic to $M_{d,r}^{\oplus e}$ where $d,r\in\BZ$, $r>0$, $(d,r)=1$, $\lambda=\frac{d}{r}$ and $re=\rk M$.}
{A $\sigma$-module $M$ over $\bigcup_r L\ancon[r]$ is called \emph{(isoclinic of) slope $\lambda$} if $M\otimes_{\bigcup_r L\ancon[r]}\bigcup_r C\ancon[r]$ is isomorphic to $M_{d,r}^{\oplus e}$ where $d,r\in\BZ$, $r>0$, $(d,r)=1$, $\lambda=\frac{d}{r}$ and $re=\rk M$.}

\subsection{The Slope Filtration Theorem} \label{SectSFT}

\comp{Let $M$ be a $\phi$-module over $\bB^\dagger_{K,\rig}$. Then there exists a unique filtration $0=M_0\subset M_1\subset\ldots\subset M_\ell=M$ of $M$ by saturated $\phi$-submodules, such that the quotients $M_i/M_{i-1}$ are isoclinic of some slopes $\lambda_i$, and $\lambda_1<\ldots<\lambda_\ell$, see Kedlaya~\cite[Theorem 6.4.1]{Kedlaya}.}
{Let $M$ be a $\sigma$-module over $\bigcup_r L\ancon[r]$. Then there exists a unique filtration $0=M_0\subset M_1\subset\ldots\subset M_\ell=M$ of $M$ by saturated $\sigma$-submodules, such that the quotients $M_i/M_{i-1}$ are isoclinic of some slopes $\lambda_i$, and $\lambda_1<\ldots<\lambda_\ell$, see Hartl~\cite[Theorem 1.7.7]{HartlPSp}.}

\comp{The base change functor from isoclinic $\phi$-modules over $\bB^\dagger_K$ of slope $\lambda$ to isoclinic $\phi$-modules over $\bB^\dagger_{K,\rig}$ is an equivalence of categories; see \cite[Theorem 6.3.3]{Kedlaya}.}
{The base change functor from isoclinic $\sigma$-modules over $\bigcup_r L\con[r]$ of slope $\lambda$ to isoclinic $\sigma$-modules over $\bigcup_r L\ancon[r]$ is an equivalence of categories; see \cite[Theorem 1.7.5]{HartlPSp}.}

\comp{Due to the limitations of the Theory of the field of norms it is \emph{not} yet possible to prove a relative version of this equivalence.}
{There \emph{exists} a relative version of this equivalence over arbitrary affinoid $L$-algebras $R$ in place of $L$. This is the main tool to prove the Theorem in \ref{SectRZConj}.}

\section{Galois Representations} \label{ChaptGaloisRep}

\comp{Let $K$ be as in Section~\ref{ChaptFontainesRings} and denote by $\CG_K:=\Gal(K^\alg/K)$ the absolute Galois group of $K$. Fontaine Theory is the theory of the category $\Rep_{\BQ_p}\CG_K$ of continuous representations of $\CG_K$ in finite dimensional $\BQ_p$-vector spaces and various full subcategories like crystalline, or semi-stable, or de Rham representations, etc.}
{Let $L$ be as in Section~\ref{ChaptFontainesRings} and let $\CG_L:=\Gal(L^\sep/L)$ be the absolute Galois group of $L$. The theory of the category $\Rep_{\BF_q\dpl z\dpr}\CG_L$ of continuous representations of $\CG_L$ in finite dimensional $\BF_q\dpl z\dpr$-vector spaces is spoiled by various unpleasant facts, such as the ones mentioned in \ref{SectAx}, \ref{SectTateTwists}, \ref{SectWTBMax}, \ref{SectBdR}, or \ref{SectOverconvergence} below.}

\subsection{$(\phi,\Gamma)$-Modules} \label{SectPhiGamma}

\comp{By the Theorem of the field of norms (see \ref{SectFieldOfNorms}) and Katz'~\cite[Proposition 4.1.1]{Katz1} the category $\Rep_{\BQ_p}\CH_K$ is equivalent to the category of slope zero $\phi$-modules over $\bB_K$ (see \ref{SectNonPerfect}) by mapping the $\CH_K$-representation $V$ to the $\phi$-module $\bD(V):=(\bB\otimes_{\BQ_p}V)^{\CH_K}$ and conversely mapping the $\phi$-module $M$ to the $\CH_K$-representation $V_p(M):=(\bB\otimes_{\bB_K}M)^{\phi=1}$.}
{By the analogue of Katz' theorem the category $\Rep_{\BF_q\dpl z\dpr}\CG_L$ is equivalent to the category of slope zero $\sigma$-modules over $L\dpl z\dpr$ by mapping the $\CG_L$-representation $V$ to the $\sigma$-module $\bigl(L^\sep\dpl z\dpr\otimes_{\BF_q\dpl z\dpr}V\bigr)^{\CG_L}$ and conversely mapping the $\sigma$-module $M$ to the $\CG_L$-representation $V_z(M):=\bigl(L^\sep\dpl z\dpr\otimes_{L\dpl z\dpr}M\bigr)^{\sigma=1}$.}

\comp{There is a relative version for certain affinoid $K$-algebras by Andreatta and Brinon~\cite{AB}.}
{There is a relative version for arbitrary affinoid $L$-algebras $R$ by replacing $\CG_L$ with the \'etale fundamental group $\pi_1^\et(\Spm R)$ (see~\cite{dJ}) and $L\dpl z\dpr$ with $R\dbl z\dbr[\frac{1}{z}]$.}

\comp{As a consequence the category $\Rep_{\BQ_p}\CG_K$ is equivalent to the category of slope zero $(\phi,\Gamma_K)$-modules over $\bB_K$, that is slope zero $\phi$-modules with $\Gamma_K:=\CG_K/\CH_K$-action commuting with $\phi$. Namely, $V$ is mapped to $\bD(V)=(\bB\otimes_{\BQ_p}V)^{\CH_K}$ which inherits the $\Gamma_K$-action from $\CG_K$.}
{Since the distinction between $\CG_L$ and $\CH_L$ collapses there is no $\Gamma_L$-action, see \ref{SectFieldOfNorms}.}

\subsection{Overconvergence} \label{SectOverconvergence}

\comp{The main theorem of Cherbonnier-Colmez~\cite{CC} says that every representation $V$ in $\Rep_{\BQ_p}\CG_K$ is \emph{overconvergent}, that is that $\bD(V)$ has a basis consisting of elements of $\bD^\dagger(V):=(\bB^\dagger\otimes_{\BQ_p}V)^{\CH_K}$, see \ref{SectOverconRings}. Equivalently every slope zero $(\phi,\Gamma_K)$-module over $\bB_K$ has a basis on which $\phi$ acts by a matrix with coefficients in $\bB^\dagger_K$.}
{Overconvergence fails in equal characteristic as the slope zero $\sigma$-module $(M,F_M)$ with 
\[
M=\BFZ\dpl z\dpr\,,\quad F_M=\bigl(\sum_{i=0}^\infty\zeta^{q^{-i}}z^i\bigr)\cdot\sigma
\]
shows. Since in mixed characteristic overconvergence is most important for the theory of $p$-adic Galois representations we propose to view $\sigma$-modules over $\bigcup_r R\con[r]$, see \ref{SectOverconRings}, as the appropriate analogue in equal characteristic of $p$-adic Galois representations.}

\subsection{Crystalline Galois Representations} \label{SectCrystReps}

\comp{A representation $V$ in $\Rep_{\BQ_p}\CG_K$ is called \emph{crystalline} if $\dim_{\BQ_p}V$ equals the dimension over $(\wt\bB_\rig)^{\CG_K}=K_0$, see \ref{SectThePeriods}, of 
\[
\bD_\cris(V)\es:=\es(\wt\bB_\rig\otimes_{\BQ_p}V)^{\CG_K}\,.
\]
The Frobenius $\phi$ on $\wt\bB_\rig$ makes $\bD_\cris(V)$ into an $F$-isocrystal over $k$.}
{Due to the unfavorable fact that 
\[
\TS\CO_C\dbl z,\frac{1}{z}\}^{\CG_L}\es=\es\CO_{L^p}\dbl z,\frac{1}{z}\}\es\supsetneq\es \ell\dpl z\dpr\,,
\]
see \ref{SectWTBMax}, the notion of crystalline Galois representations is problematic. But see the next item.}

\subsection{Rigidified Local Shtuka} \label{SectRigLocalShtuka}

\comp{The Breuil-Kisin classification \cite[Theorem 0.1]{Kisin} of crystalline $\CG_K$-representations states the following:\\
Let $\FS:=W(k)\dbl u\dbr$ and let $\pi\in K$ be a uniformizer with Eisenstein polynomial $E(u)\in\FS$. Let $\phi:\FS\to\FS$ extend the Frobenius lift $\phi$ on $W(k)$ and map $u\mapsto u^p$. Let $\Mod^\phi_{/\FS}$ be the Tannakian category of finite free $\FS$-modules $\FM$ with an isomorphism $F_\FM:\phi^\ast\FM[\frac{1}{E(u)}]\isoto\FM[\frac{1}{E(u)}]$, where $\phi^\ast\FM:=\FM\otimes_{\FS,\phi}\FS$.}
{ A \emph{rigidified local shtuka} $(M,F_M,\delta_M)$ over $\CO_L$ consists of a local shtuka $(M,F_M)$ over $\CO_L$ and, posing $(D,F_D):=(M,F_M)\otimes_{\CO_L\dbl z\dbr}\ell\dpl z\dpr$, an isomorphism $\delta_M$  
\[
\TS M\otimes_{\CO_L\dbl z\dbr}\CO_L\dbl z,\frac{1}{z}\}[\frac{1}{t_{\SSC -}}]\isoto D\otimes_{\ell\dpl z\dpr}\CO_L\dbl z,\frac{1}{z}\}[\frac{1}{t_{\SSC -}}]
\]
which satisfies $\delta_M\circ F_M=F_D\circ\sigma^\ast\delta_M$ and which reduces to the identity modulo 
$\Fm_L$
. 
If $L$ is discretely valued the forgetful functor $(M,F_M,\delta_M)\mapsto (M,F_M)$ is an equivalence of categories by \cite[Lemma 2.1]{GL}, see also \cite[Lemma 2.3.1]{HartlPSp}.}

\comp{Then the category of crystalline representations is tensor equivalent to a full subcategory of the isogeny category of $\Mod^\phi_{/\FS}$. The essential image can be characterized in terms of a connection.}
{From a (rigidified) local shtuka $(M,F_M)$ over $\CO_L$ one obtains a Galois representation
\[
V_z(M):=\bigl(L^\sep\dpl z\dpr\otimes_{\CO_L\dbl z\dbr}M\bigr)^{\sigma=1}
\]
as in \ref{SectTateModules}.
The tensor functor $(M,F_M)\mapsto V_z(M)$ is faithful; see \cite[Proposition 2.1.4]{HartlPSp}.}

\subsection{The Mysterious Functor} \label{SectMysteriousFunctor}

\comp{If $V$ is a crystalline $\CG_K$-representation, then $\bD_\cris(V)$ is an $F$-isocrystal over $k$, see \ref{SectCrystReps}. The embedding $\wt\bB_\rig\subset\bB_\dR$ (see \ref{SectThePeriods} and \ref{SectBdR}) equips $\bD_\cris(V)\otimes_{K_0}K$ with a filtration $Fil^\bullet$ by $K$-subspaces, such that $\bigl(\bD_\cris(V),Fil^\bullet\bigr)$ is a filtered isocrystal over $K$, see \ref{SectFilteredIsocrystals}. This is Fontaine's functor
\[
V\es\mapsto\es\bigl(\bD_\cris(V),Fil^\bullet\bigr)
\]
whose existence was conjectured by Grothen\-dieck (the ``mysterious functor'') in case $V$ equals the \'etale cohomology of a smooth proper $K$-variety with good reduction.}
{If $(M,F_M,\delta_M)$ is a rigidified local shtuka over $\CO_L$ then $(D,F_D):=(M,F_M)\otimes_{\CO_L\dbl z\dbr}\ell\dpl z\dpr$ is a $z$-isocrystal over $\ell$ and 
\[
\Fq_D\es:=\es\sigma^\ast\delta_M\circ F_M^{-1}\bigl(M\otimes_{\CO_L\dbl z\dbr}L\dbl z-\zeta\dbr\bigr)
\]
 is a Hodge-Pink structure on $(D,F_D)$, such that $(D,F_D,\Fq_D)$ is a filtered isocrystal over $L$, see \ref{SectFilteredIsocrystals}.
The functor 
\[
\BH:\es(M,F_M,\delta_M)\es\mapsto\es(D,F_D,\Fq_D)
\]
is the analogue of the mysterious functor.}

\comp{A filtered isocrystal in the essential image of Fontaine's functor is called \emph{admissible}.}
{A filtered isocrystal in the essential image of $\BH$ is called \emph{admissible}.}

\subsection{Weakly Admissible Implies Admissible} \label{SectWA=>A}

\comp{By the Colmez-Fontaine Theorem~\cite{CF} Fontaine's Functor is an equivalence of categories between crystalline $\CG_K$-representations and weakly admissible filtered isocrystals over $K$, see \ref{SectWA}. (Here $K$ has to be discretely valued with perfect residue field.)}
{The functor $\BH$ is an equivalence between rigidified local shtuka over $\CO_L$ and weakly admissible filtered isocrystals (see \ref{SectWA}) if the completion $\wt L$ of the compositum $\ell^\alg L$ inside $C$ does not contain an element $a$ with $0<|a|<1$ such that for all $n\in\BN$ the $q^n$-th roots of $a$ also lie in $\wt L$, see \cite[Theorem 2.5.3]{HartlPSp}. This condition is for example satisfied if the value group of $L$ is not $q$-divisible.}

\section{Period Spaces for Filtered Isocrystals} \label{ChaptPeriodSpaces}
\subsection{Period Spaces}

\comp{On a fixed $F$-isocrystal $(D,F_D)$ over $\BF_p^{\,\alg}$ the filtrations from \ref{SectFilteredIsocrystals} are parametrized by a projective partial flag variety over $W(\BF_p^{\,\alg})[\frac{1}{p}]$.}
{On a fixed $z$-isocrystal $(D,F_D)$ over $\BF_q^{\,\alg}$ the Hodge-Pink structures from \ref{SectFilteredIsocrystals} are parametrized by a quasi-projective partial jet bundle over a partial flag variety over $\BF_p^{\,\alg}\dpl\zeta\dpr$. The jets arise since a Hodge-Pink structure contains more information than just the Hodge-Pink filtration.}

\comp{The weakly admissible filtrations (see \ref{SectWA}) form a rigid ana\-lytic subspace $\CF^{wa}$ of this partial flag variety by Rapoport-Zink~\cite[Proposition 1.36]{RZ}. This is a \emph{$p$-adic period space}. When viewed as a Berkovich space it is even a connected open Berkovich subspace of the partial flag variety; see \cite[Proposition 1.3]{HartlRZ}.}
{The weakly admissible Hodge-Pink structures (see \ref{SectWA}) form a rigid analytic subspace $\CF^{wa}$ of this partial jet bundle by \cite[Theorem 3.2.5]{HartlPSp} called a \emph{period space for Hodge-Pink structures}. When viewed as a Berkovich space it is even a connected open Berkovich subspace of the partial jet bundle.}

\subsection{A Conjecture of Rapoport and Zink} \label{SectRZConj}

\comp{Rapoport and Zink~\cite[p.\ 29]{RZ} conjecture the existence of an \'etale morphism $\CF'\to\CF^{wa}$ of rigid analytic spaces, which is bijective on $K$-valued points with $K$ finite over $W(\BF_p^{\,\alg})[\frac{1}{p}]$, and a $p$-adic representation of the \'etale fundamental group
\[
\pi_1^\et(\CF')\es\longto\es\GL_n(\BQ_p)
\]
which induces the universal filtered isocrystal over $\CF'$.}
{The analogue of Rapoport and Zink's conjecture is a theorem; see \cite[Theorem 3.4.3]{HartlPSp}:\\
There exists a unique maximal open Berkovich subspace $\CF^a$ of the period space $\CF^{wa}$ containing all its $L$-valued points for fields $L$ as in \ref{SectWA=>A}, there exists an admissible formal scheme $Y$ over $\BF_q^{\,\alg}\dbl\zeta\dbr$ whose associated Berkovich space $Y^{\rm Berko}$ is an \'etale covering space of $\CF^a$, and there exists a rigidified local shtuka (see \ref{SectRigLocalShtuka}) over $Y$ inducing the universal filtered isocrystal over $\CF^a$. This local shtuka gives rise to a representation $\pi_1^\et(\CF^a)\to\GL_n\bigl(\BF_q\dpl z\dpr\bigr)$ of the \'etale fundamental group. The proof relies on the relative descent mentioned in \ref{SectSFT} which is lacking in mixed characteristic.}

\subsection{Rapoport-Zink Spaces} \label{SectRZSpaces}

\comp{Let $\BE$ be a fixed $p$-divisible group over $\BF_q^{\,\alg}$. The functor which assigns to every formal $W(\BF_p^{\,\alg})$-scheme $S$ the set of isomorphism classes of pairs $(E,\rho)$ where
\begin{itemize}
\item 
$E$ is a $p$-divisible group over $S$,
\item 
$\rho:\BE_{\bar S}\to E_{\bar S}$ is a quasi-isogeny over $\bar S$, the closed subscheme of $S$ defined by the ideal $p\,\CO_S$,
\end{itemize}
is representable by an adic formal scheme $X$ locally formally of finite type; see \cite[Theorem 2.16]{RZ}. $X$ is called a \emph{Rapoport-Zink space}.}
{Let $\BM$ be a fixed local shtuka over $\BF_q^{\,\alg}$. The functor which assigns to every formal $\BF_q^{\,\alg}\dbl\zeta\dbr$-scheme $S$ the set of isomorphism classes of pairs $(M,\rho)$ where
\begin{itemize}
\item 
$M$ is a local shtuka over $S$,
\item 
$\rho:M\otimes_{\CO_S\dbl z\dbr}\CO_{\bar S}\dbl z\dbr[\frac{1}{z}]\isoto\BM_{\BF_q^{\,\alg}\dbl z\dbr}\CO_{\bar S}\dbl z\dbr[\frac{1}{z}]$ is an isomorphism where $\bar S$ is the closed subscheme of $S$ defined by the ideal $\zeta\,\CO_S$,
\end{itemize}
is representable by an adic formal scheme $X$ locally formally of finite type; see \cite
{Crystals}. $X$ is called a \emph{Rapoport-Zink space}.}

\subsection{Period Morphisms} \label{SectPeriodMorph}

\comp{The Hodge-Tate filtration on the Dieudonn\'e crystal associated with the universal $p$-divisible group over the formal scheme $X$ from \ref{SectRZSpaces} defines an \'etale \emph{period morphism} $X^\rig\to\CF^{wa}$ by \cite[Proposition 5.15]{RZ}. When viewed as a morphism of Berkovich spaces it has open image and identifies the Berkovich space associated with $X^\rig$ with an \'etale covering space of this image.}
{Applying the mysterious functor from \ref{SectMysteriousFunctor} to the universal local shtuka over the formal scheme $X$ from \ref{SectRZSpaces} defines an \'etale \emph{period morphism} from $X^{\rm Berko}$ to the space $\CF^a$ from \ref{SectRZConj}. It identifies $X^{\rm Berko}$ with the \'etale covering space $Y^{\rm Berko}$ of $\CF^a$ from \ref{SectRZConj}, see \cite{Crystals}.}

\vspace{2cm}


\end{center}

%
%

\vfill

\noindent
Urs Hartl\\
University of Freiburg\\
Institute of Mathematics\\
Eckerstr.~1\\
D -- 79104 Freiburg\\
Germany\\[1mm]
E-mail: 
urs.hartl@math.uni-freiburg.de


\begin{thebibliography}{EGA}
\addcontentsline{toc}{section}{References}

 
\bibitem{Anderson} G.\ Anderson: $t$-Motives, {\it Duke Math.\ J.\/} {\bf 53} (1986), 457--502.


\bibitem{Anderson2} G.\ Anderson: On Tate Modules of Formal $t$-Modules, {\em Internat.\ Math.\ Res.\ Notices} {\bf 2} (1993), 41--52.

\bibitem{Andreatta} F.\ Andreatta: {\em Generalized ring of norms and generalized $(\phi,\Gamma)$-modules}, University of Padova, Preprint 2005.

\bibitem{AB} F.\ Andreatta, O.\ Brinon: {\em Surconvergence des r\'epresentations $p$-adiques: le cas relativ}, Preprint April 2006.

\bibitem{Ax} J.\ Ax: Zero of polynomials over local fields -- the Galois action, {\em J.\ Algebra} {\bf 15} (1970), 417--428.



 


\bibitem{Bosch} S.\ Bosch: {\em Lectures on Formal and Rigid Geometry}, Preprint {\bf 378}, University of M\"unster, SFB 478--Preprint Series, M\"unster 2005, available at \\
{\tt http://wwwmath1.uni-muenster.de/sfb/about/publ/heft378.ps}\,.

\bibitem{BGR} S.\ Bosch, U.\ G{\"u}ntzer, R.\ Remmert: {\em Non-archimedean analysis}, Grundlehren {\bf 261}, Springer-Verlag, Berlin-Heidelberg 1984. 
 
\bibitem{FRG} S.\ Bosch, W.\ L{\"u}tkebohmert: Formal and Rigid Geometry  
I.\ Rigid Spaces, {\em Math.\ Ann.\/} {\bf 295} (1993), 291--317. 




\bibitem{Brinon} O.\ Brinon: {\em Repr\'esentations Galoisiennes $p$-adiques dans le cas relative}, Th\`ese de l'Universit\'e de Paris 11 (Orsay), 2004.

\bibitem{CC} F.\ Cherbonnier, P.\ Colmez: Repr\'esentations $p$-adiques surconvergentes, {\em Invent.\ Math.} {\bf 133} (1998), no.\ 3, 581--611.




\bibitem{Colmez02} P.\ Colmez: Espaces de Banach de dimension finie, {\em J.\ Inst.\ Math.\ Jussieu} {\bf 1} (2002), no.\ 3, 331--439.

\bibitem{Colmez} P.\ Colmez: {\em Espaces Vectoriels de dimension finie et r\'epresentations de de Rham}, Preprint 2003, available at {\tt http://www.math.jussieu.fr/$\sim$colmez}.

\bibitem{CF} P.\ Colmez, J.-M.\ Fontaine: Construction des repr\'esentations $p$-adiques semi-stables, {\em Invent.\ Math.} {\bf 140} (2000), no. 1, 1--43.



\bibitem{Demazure} M.\ Demazure: {\em $p$-divisible groups}, LNM {\bf 302}, Springer-Verlag, Berlin-New York, 1972.

\bibitem{Drinfeld} V.G.\ Drinfeld: Elliptic Modules, {\it Math.\ USSR-Sb.\/} {\bf 23} (1976), 561--592.

 


\bibitem{Faltings1} G.\ Faltings: $p$-adic Hodge theory, {\em J.\ Amer.\ Math.\ Soc.} {\bf 1}  (1988),  no.\ 1, 255--299.

\bibitem{Faltings} G.\ Faltings: Crystalline cohomology and $p$-adic Galois-representations, {\em Algebraic analysis, geometry, and number theory (Baltimore 1988)}, pp.\ 25--80, Johns Hopkins Univ.\ Press, Baltimore 1989.


\bibitem{Fontaine77} J.-M.\ Fontaine: {\em Groupes $p$-divisibles sur les corps locaux} Ast\'erisque {\bf 47-48} (1977).






\bibitem{GL} A.\ Genestier, V.\ Lafforgue: Th\'eories de Fontaine et Fontaine-Laffaille en \'egale charact\'eristique, Preprint 2004.



\bibitem{Grothendieck} A.\ Grothendieck: \emph{Groupes de Barsotti-Tate et cristaux de Dieudonn\'e}, S\'eminaire de Math\'ematiques Sup\'erieures 45, Presses de l'Universit\'e de Montr\'eal, Montr\'eal 1974.

\bibitem{GossDict} D.\ Goss: Dictionary, {\em The arithmetic of function fields (Columbus, OH, 1991)},  pp.\ 475--482, Ohio State Univ.\ Math.\ Res.\ Inst.\ Publ., 2, de Gruyter, Berlin, 1992.

\bibitem{Harder} G.\ Harder: Wittvektoren, {\em Jahresber.\ Deutsch.\ Math.-Verein.} {\bf 99} (1997), no.\ 1, 18--48.

\bibitem{AbSh} U.\ Hartl: Uniformizing the Stacks of Abelian Sheaves, in \emph{Number Fields and Function fields - Two Parallel Worlds, Papers from the 4th Conference held on Texel Island, April 2004}, G.\ van der Geer, B.\ Moonen, R.\ Schoof, Editors, pp.\ 167--222, Progress in Math.\ 239, Birkh\"auser-Verlag, Basel 2005. See also {\tt arXiv:math.NT/0409341}. 

\bibitem{HartlPSp} U.\ Hartl: {\em Period Spaces for Hodge Structures in Equal Characteristic}, Preprint on {\tt arXiv:math.NT/0511686}.

\bibitem{HartlRZ} U.\ Hartl: {\em On a Conjecture of Rapoport and Zink}, Preprint on {\tt arXiv:math.NT/0605254}.

\bibitem{Crystals} U.\ Hartl: {\em Local Shtuka and Divisible Local Anderson Modules}, in preparation.

\bibitem{HP} U.\ Hartl, R.\ Pink: Vector bundles with a Frobenius structure on the punctured unit disc, {\em Comp.\ Math.} {\bf 140} n.3 (2004), 689--716. 

 

\bibitem{Illusie} L.\ Illusie: Cohomologie de de Rham et cohomologie \'etale $p$-adique (d'apr\`es G.\ Faltings, J.-M.\ Fontaine et al.), {\em S\'eminaire Bourbaki, Vol. 1989/90,  Ast\'erisque} {\bf 189-190} (1990), Exp.\ No.\ 726, 325--374.

\bibitem{dJ} J.\ de Jong: {\'E}tale Fundamental groups of non-Archimedean analytic spaces, {\em Comp.\ Math.\/} {\bf 97} (1995), 89--118. 


 

\bibitem{Katz1} N.\ Katz: $p$-adic properties of modular schemes and modular forms, {\em Modular functions of one variable, III (Proc.\ Internat.\ Summer School, Univ.\ Antwerp, Antwerp, 1972)},  pp.\ 69--190. LNM {\bf 350}, Springer, Berlin, 1973. 

\bibitem{Katz} N.\ Katz: Slope filtration of $F$-crystals.  {\em Journ{\'e}es de G{\'e}om{\'e}trie Alg{\'e}brique de Rennes (Rennes, 1978)}, Vol. I,  pp. 113--163, Ast{\'e}risque {\bf 63}, Soc.\ Math.\ France, Paris 1979.

 

\bibitem{Kedlaya} K.\ Kedlaya:  Slope filtrations revisited, {\em Doc.\ Math.} {\bf 10} (2005), 447--525.

\bibitem{Kisin} M.\ Kisin: {\em Crystalline Representations and $F$-Crystals}, Preprint on {\tt http://math.uchicago.edu/$\sim$kisin}, version from April 8th, 2005. 



\bibitem{Laumon} G.\ Laumon: {\it Cohomology of Drinfeld Modular Varieties I\/}, Cambridge Studies in Advanced Mathematics {\bf 41}, Cambridge University Press, Cambridge, 1996.

\bibitem{Laumon97} G.\ Laumon: Drinfeld shtukas, in {\em Vector bundles on curves---new directions (Cetraro, 1995)}, pp.\ 50--109, LNM {\bf 1649}, Springer, Berlin, 1997.


\bibitem{LT} J.\ Lubin, J.\ Tate: Formal complex multiplication in local fields, {\em Ann.\ of Math.} (2) {\bf 81} (1965), 380--387.

 
\bibitem{Manin} Y.I.\ Manin: The theory of commutative formal groups over fields of finite characteristic, {\it Russian Math. Surveys} {\bf 18} No.~6 (1963), 1--83.

\bibitem{Messing} W.\ Messing: {\em The crystals associated to Barsotti-Tate groups: with applications to abelian schemes}, LNM {\bf 264}, Springer-Verlag, Berlin-New York, 1972.

\bibitem{Pink} R.\ Pink: {\em Hodge Structures over Function Fields}, Preprint 1997, available at {\tt http://www.math.ethz.ch/$\sim$pink}.



\bibitem{RZ} M.\ Rapoport, T.\ Zink: {\it Period Spaces for $p$-divisible Groups\/}, Ann.\ Math.\ Stud.\ {\bf 141}, Princeton University Press, Princeton 1996.

 



\bibitem{Serre} J.-P.\ Serre: {\em Local fields}, GTM {\bf 67}, Springer-Verlag, New York-Berlin, 1979.

\bibitem{Tate2} J.\ Tate: $p$-divisible groups, {\em Proc.\ Conf.\ Local Fields (Driebergen, 1966)}, pp. 158--183, Springer-Verlag, Berlin 1967.


\bibitem{Tsuji} T.\ Tsuji: $p$-adic \'etale cohomology and crystalline cohomology in the semi-stable reduction case, {\em Invent.\ Math.} {\bf 137}  (1999),  no.\ 2, 233--411.


\bibitem{Wintenberger94} J.-P.\ Wintenberger: Th\'eor\`eme de comparaison $p$-adique pour les sch\'emas ab\'eliens, I. Construction de l'accouplement de p\'eriodes, {\em P\'eriodes $p$-adiques (Bures-sur-Yvette, 1988).  Ast\'erisque} {\bf 223} (1994), 349--397.


\end{thebibliography}
\end{document}